\documentclass{article}
\usepackage{amsmath,amsthm,amsfonts,amscd,amssymb,eucal,latexsym}

\newcommand\cA{{\cal A}}
\newcommand\cB{{\cal B}}
\newcommand\cC{{\cal C}}

\newcommand\cH{{\cal H}}
\newcommand\cI{{\cal I}}
\newcommand\cJ{{\cal J}}

\newcommand\cM{\cal {M}}

\newcommand\cQ{{\cal Q}}

\newcommand\cS{{\cal S}}
\newcommand\cT{{\cal T}}

\begin{document}

\title{Ergodic Actions of Compact Quantum Groups from Solutions of the Conjugate Equations} 

\author{Claudia Pinzari$^a$, John E. Roberts$^b$}
\normalsize
\author{Claudia Pinzari$^a$, John E. Roberts$^b$\\ \\
$^a\,$Dipartimento di Matematica, Universit\`a di Roma ``La Sapienza''\\
00185--Roma, Italy\\
$^b\,$Dipartimento di Matematica, Universit\`a di Roma ``Tor Vergata'',\\
00133 Roma, Italy}
\date{}
\maketitle

\begin{abstract} \noindent
We use a  tensor $C^*$--category with conjugates and two quasitensor functors 
into the category of Hilbert spaces to define a ${}^*$--algebra depending 
functorially on this data. In one of them is tensorial we can complete in the maximal $C^*$--norm. A particular case of this construction allows us 
to begin with solutions of the conjugate equations and associate ergodic actions 
of quantum groups on the $C^*$--algebra in question. The quantum groups 
involved are $A_u(Q)$ and $B_u(Q)$.
\begin{footnote}{2000 American Mathematical Society Classification: Primary 18D10, 
20N99, 37A55}\end{footnote} 
\begin{footnote}{Keywords: ergodic actions, quantum groups, conjugate equations}
\end{footnote}
 
\end{abstract}

\begin{section} {Introduction}
The theory of ergodic actions of compact quantum groups on 
unital $C^*$--algebras
has recently attracted interest. In the group case,  one of the first 
results was the theorem by H\o egh-Krohn, Landstad and St\o rmer 
 asserting 
that the multiplicity of an irreducible representation 
is always bounded by its dimension and that the unique $G$--invariant 
state is a trace \cite{HLS}.

Ergodic theory for group actions was later investigated by 
Wassermann in a series of papers \cite{Wassermann1}, 
\cite{Wassermann2}, \cite{Wassermann3}, who, among other 
results,  classified all ergodic actions of $SU(2)$ on 
von Neumann algebras. In particular, he proved the 
important result
that $SU(2)$ cannot act ergodically on the hyperfinite 
$II_1$ factor.

For compact quantum groups, ergodic theory on $C^*$--algebras was 
initiated by Boca. He generalized the HLS theorem, showing that 
the multiplicity of an irreducible 
is bounded instead by its quantum dimension. Woronowicz
\cite{Wcmp} noticed that the modular group of a compact quantum group 
is not always trivial, consequently, the invariant state cannot be a trace in general. 
Boca  described the  modularity of this state for a general 
ergodic action \cite{Boca}.

Podles was interested in studying quantum spheres: he  
introduced subgroups and quotients for compact quantum groups 
and computed them for the quantum $SU(2)$ and $SO(3)$ groups. Some 
of the quantum spheres he found are not embedable into the quotient 
spaces \cite{Podles}.

Later Wang  found many  examples of  ergodic quantum actions 
on $C^*$--algebras: 
for the quantum groups $A_u(Q)$ on type $III_\lambda$ 
Powers factors,
on the Cuntz algebras,   on the 
injective factor of type $III_1$ and on the hyperfinite $II_1$ factor 
(this, by a Kac type quantum group). He also
found an example 
on a commutative 
$C^*$--algebra that is not a quotient \cite{Wang}.

Classifying ergodic $C^*$--actions of the quantum group 
$S_\mu U(2)$ of  Woronowicz is an 
open problem. Tomatsu has classified all 
those which are embedded in the translation action of
$S_\mu U(2)$ \cite{Tomatsu}.

Bichon, De Rijdt and Vaes have constructed examples of ergodic actions of 
 $S_\mu U(2)$ not embedable in the translation action. 
since the multiplicity of an irreducible is bigger
than its integral dimension. 
The 
authors also introduced a new 
invariant, the quantum multiplicity $m(u)$ of an irreducible 
representation $u$ in 
the action. 
This invariant reduces to the quantum dimension for the translation action.
In general, one has the bounds: $\text{multiplicity}(u)\leq m(u)\leq 
\text{q-dim}(u)$.
Even for quotient actions, the quantum multiplicity is not 
an integer in general \cite{BRV}. These examples were constructed by means of a generalisation of the Tannaka-Krein duality theorem for compact quantum groups \cite{WoronowiczTK} to ergodic actions of full multiplicity.

In \cite{PR} we in turn extended the duality theorem of \cite{BRV} to general ergodic actions of compact quantum groups on unital $C^*$--algebras. To this aim, we introduced the notion of quasitensor functor between two tensor $C^*$--categories with conjugates, and we showed that quasitensor functors from the representation category of a compact quantum group $G$ to the category of Hilbert spaces characterise the spectral functors of ergodic $C^*$--actions of $G$.

As an application, we constructed ergodic actions of  $S_\mu U(d)$  starting from abstract 
tensor $C^*$--categories with a Hecke symmetry of parameter $q=\mu^2$.
In particular, one gets ergodic actions of $S_\mu 
U(2)$ from a real or pseudoreal object $y$ of a tensor $C^*$--category 
with 
intrinsic dimension  $d(y)\geq2$, with
 $\mu$ and $d(y)$
 related by 
$d(y)=|\mu+\mu^{-1}|$ and $\mu$ positive if  $y$ is 
pseudoreal and 
negative otherwise. 

Our interest in braiding was motivated by 
low 
dimensional QFT, where
braided tensor $C^*$--categories arise, albeit with a {\it unitary} 
braiding \cite{Haag}.

The aim of this paper is twofold. We first give an alternative notion of quasitensor functor  and we show the equivalence with that of \cite{PR}, as well as the construction of the   mentioned $C^*$--ergodic action of $G$. 
Furthermore, we apply this construction to obtain   ergodic actions of  compact quantum groups starting from solutions of 
the conjugate equation in a tensor $C^*$--category.

This paper has a sequel, that is \cite{PRinduction},    where we start  precisely from the non commutative space  here  obtained  to construct Hilbert bimodule representations of compact quantum groups arising from tensor $C^*$--categories
generated by an object of intrinsic dimension $\geq 2$. In this sense,
  our ergodic actions should be regarded as virtual quantum subgroups.

Furthermore, as an application, we get ergodic actions of $S_{\mu}U(2)$, for  {negative} or positive values
of the deformation parameter uniquely determined by the intrinsic dimension.
Ergodic actions of the quantum  groups $B_u(Q)$ and $A_u(Q)$ of Wang appear as well.

The paper is organized as follows. In Section 2 we review the main facts 
about compact quantum groups \cite{WLesHouches}, the main invariants of 
ergodic $C^*$--actions and the duality theorem of \cite{PR}.

In Section 3 we state our main results: the existence of ergodic 
$C^*$--actions 
of the compact quantum group $A_u(Q)$, in the notation of 
Wang, associated with an invertible positive matrix $Q\in M_n({\mathbb C})$ 
with Tr$(Q)=Tr(Q^{-1})$ arising from normalized solutions $R,\overline R$ of the conjugate equations.
$Q$ and $R$ are related by Trace$(Q)=R^*\circ R$.  

For self-conjugate solutions of the conjugate equations we also get 
ergodic $C^*$--actions of the compact quantum group $B_u(Q)$, in the  
notation of Wang, associated with an invertible matrix $Q\in M_n(\mathbb C)$ 
with  $Q\overline Q=\pm 1$ and hence Tr$Q^*Q=$Tr$(Q^*Q)^{-1}$. $Q$ and  $R$ 
are now related by Trace$(Q^*Q)=R^*\circ R$ and $\pm 1$ 
distinguishes real from pseudoreal solutions.
\begin{footnote} {Banica denotes $A_u(Q)$ by $A_u(F)$ where $Q=F^*F$ and 
$B_u(Q)$ by $A_o(F)$, where $Q=F^*$, (D\' efinition 2 of \cite{Banica}).}\end{footnote}

 In Section 4, we generalize, at an algebraic level, the construction of the ${}^*$--algebra 
 carrying the ergodic action giving a formalism symmetric in two 
 quasitensor functors and discuss its functorial properties. We then complete in the maximal $C^*$--norm when one of the two functors is tensorial. The general case will be considered elsewhere \cite{PRlinking}. The 
 action itself is defined in Section 5. Section 6 recalls some properties 
 of the Temperley-Lieb categories associated with self-conjugate 
 solutions of the conjugate equations whilst Section 7 is devoted 
 to related categories needed for treating general solutions of the 
 conjugate equations. Section 8 treats the embeddings of these 
 categories into the category of Hilbert spaces and the associated 
 compact quantum groups and concludes with the proof of 
 the main results.\medskip

 \end{section}

\begin{section}{Preliminaries}

In this preliminary section we  recall the  main invariants 
associated with an ergodic action of a compact quantum group $G$ on a 
unital $C^*$--algebra ${\cal C}$ and the duality theorem of \cite{PR}.
\medskip

\noindent{\bf 2.1 Compact quantum groups}.
\medskip

We follow  Woronowicz \cite{WLesHouches} in defining a 
compact quantum group $G$ to be a 
pair $G=({\cal Q},\Delta)$ where ${\cal Q}$ is a $C^*$--algebra with 
unit $I$
and $\Delta:{\cal Q}\to{\cal Q}\otimes{\cal Q}$ a unital coassociative
$^*$--homomorphism, {\it the coproduct}:
$$\Delta\otimes\iota\circ\Delta=\iota\otimes\Delta\circ\Delta,$$
with $\iota:{\cal Q}\to{\cal Q}$ the identity map. 
To economize on brackets we shall always evaluate tensor products 
before composition.
The coproduct 
is required  to  
satisfy the following  nondegeneracy condition: the subspaces
$I\otimes{\cal 
Q}\Delta({\cal Q})$ and ${\cal Q}\otimes I\Delta({\cal Q})$ are 
dense   in ${\cal Q}\otimes{\cal Q}$.

A unitary {\it representation} of $G$ on a finite dimensional Hilbert space 
$H_u$ is a linear map 
$u: H_u\to H_u\otimes{\cal Q}$
satisfying the group homomorphism property, nondegeneracy and unitarity, 
expressed respectively by:
$$u\otimes\iota\circ u=\iota\otimes\Delta\circ u,$$
$$u(H_u)({I\otimes\cal Q})=H_u\otimes{\cal Q}$$
$$(u(\psi), u(\psi'))_{\cal Q}=(\psi,\psi')I,$$
where on the left hand side of the last relation we have the natural 
${\cal Q}$--valued inner
product of the right Hilbert module $H_u\otimes{\cal Q}$:
$$(\psi\otimes q, \psi'\otimes q')_{\cal Q}:=(\psi,\psi')q^*q',\quad 
\psi,\psi'\in H_u, q,q'\in{\cal Q}.$$
The {\it coefficients} of $u$ are elements of ${\cal Q}$ defined by 
$u_{\psi,\psi'}:= \psi^*\otimes Iu(\psi')$, where $\psi^*: H_u\to{\mathbb 
C}$
is the annihilation operator $\psi^*\psi':=(\psi,\psi')$.
Representation coefficients span a dense $^*$--subalgebra of ${\cal Q}$.
If $u$ and $v$ are two representations, we can form the 
{\it tensor 
product} representation $u\otimes v$ on the tensor product Hilbert space
$H_u\otimes H_v$, defined by
$$(u\otimes v)_{\psi\otimes\phi, 
\psi'\otimes\phi'}:=u_{\psi,\psi'}v_{\phi,\phi'}.\eqno(2.1)$$
A {\it conjugate} of $u$ is a unitary representation $\overline{u}$ with
an antilinear invertible $j:H_u\to H_{\overline{u}}$ such that 
$$\overline{u}_{\phi,j\psi}=(u_{j^*\phi,\psi})^*.\eqno(2.2)$$
A conjugate $\overline{u}$ of $u$ is defined up to unitary equivalence.
Every representation has a conjugate representation \cite{Wcmp}.
The category $\text{Rep}(G)$ with objects unitary representations of $G$ 
and arrows the intertwining operators
$$(u,v):=\{A:H_u\to H_v: v\circ A=A\otimes I\circ u\},$$
is  a tensor $C^*$--category with conjugates (and also subobjects and 
direct 
sums) in the sense of \cite{LR}.
Furthermore,
$\text{Rep}(G)$ embeds naturally as
a tensor $^*$--subcategory of the category ${\cal H}$ of Hilbert 
spaces.
Conversely, any tensor $^*$--subcategory of ${\cal H}$ with conjugation, 
subobjects and direct sums is the representation category of a compact 
quantum group \cite{WoronowiczTK}.

The quantum groups of interest here are the Woronowicz 
deformations $S_\mu U(2)$ by a nonzero real parameter $\mu$ 
\cite{WoronowiczTK} and the Van Daele-Wang orthogonal groups and 
unitary groups $B_u(Q)$ and $A_u(Q)$ associated with
an invertible matrix $Q\in M_n({\mathbb C})$ \cite{WVD}.
\medskip

\noindent{\bf 2.2 Ergodic $C^*$--actions}
\medskip
 
Consider a unital $C^*$--algebra ${\cal C}$ and a compact quantum group 
$G$. An {\it action} of $G$ on ${\cal C}$ is a unital $^*$--homomorphism
$\delta:{\cal C}\to{\cal C}\otimes{\cal Q}$ satisfying the group 
representation property:
$$\iota\otimes\Delta\circ\delta=\delta\otimes\iota\circ\delta,$$
and the nondegeneracy property requiring that $\delta({\cal 
C})I\otimes{\cal 
Q}$ be dense in ${\cal C}\otimes{\cal Q}$. The {\it spectrum} of $\delta$,
$\text{sp}({\delta})$, is defined  to be the set of all unitary
representations $u$ of $G$ for which there is a faithful linear map
$T: H_u\to{\cal C}$ intertwining the representation $u$ with the action 
$\delta$:
$$\delta\circ T=T\otimes\iota\circ u.$$
In other words, if $u_{ij}$ are the 
coefficients of $u$ in some orthonormal basis of $H$,
we are requiring the existence of a spectral multiplet of  linearly independent elements
$c_1,\dots, c_d\in{\cal C}$, with $d$ the dimension of $u$, transforming 
like $u$ under the action:
$\delta(c_i):=\sum_j c_j\otimes u_{ji}$. The linear span of all the $c_i$'s, 
denoted ${\cal C}_{\text{sp}}$,
as $u$ varies in the spectrum, is a dense $^*$--subalgebra of ${\cal C}$
\cite{Podles}.
 
The action $\delta$ is called {\it ergodic} if the fixed point algebra
$${\cal C}^\delta=\{c\in{\cal C}:\delta(c)=c\otimes I\}$$ reduces to the 
complex numbers: ${\cal C}^\delta={\mathbb C}I.$
The simplest example of an ergodic action is the translation action of $G$ 
on
${\cal C}={\cal Q}$ with $\delta=\Delta$. Another simple class of examples 
are the adjoint actions on $B(H_u)$, where $u$ is an irreducible unitary 
representation. The spectrum then consists of the subrepresentations of 
$u\otimes\overline u$.\medskip

If an action $\delta $ is ergodic,   
the spectral 
multiplets transforming like $u$ form 
Hilbert spaces.
In fact, for any representation $u$, consider the space
$$L_u:=\{T:H_u\to{\cal C}, \delta\circ T =T\otimes\iota\circ u\}.$$
If $S,T\in L_u$, 
$<S,T>:=\sum_iT(\psi_i)S(\psi_i)^*$, where $(\psi_i)$ is an 
orthonormal basis 
of $H_u$, is an element of the 
fixed point 
algebra
${\cal C}^\delta$, and hence a complex number.
$L_u$ is known to be finite dimensional and is therefore a 
Hilbert space with the above inner product.
This Hilbert space is nonzero precisely when $u$ contains a 
subrepresentation $v\in\text{sp}(\delta)$. In particular, for an 
irreducible $u$, the conditions $u\in\text{sp}(\delta)$ and $L_u\neq 0$ 
are equivalent. The dimension of $L_u$ is called the {\it multiplicity}
of $u$ and denoted $\text{mult}(u)$.

The complex conjugate vector space $\overline{L_u}$,
 endowed with the conjugate inner product 
$$<\overline{S},\overline{T}>:=<T,S>=\sum_i 
S(\psi_i)T(\psi_i)^*,$$
is called the {\it spectral space} associated with $u$.

If for example  $\delta$ is the translation action on ${\cal Q}$, any
$\psi\in H_u$ defines an element of $L_u$ by
$$T_\psi(\psi'):=\psi^*\otimes Iu(\psi').$$
Hence the spectral space $\overline{L_u}$ 
can be identified  with $H_u$ through the  unitary map
$$\psi\in H_u\to \overline{T_\psi}\in \overline{L_u}.$$
For a general 
ergodic action, we introduce
certain maps whose coefficients generate the dense 
$^*$--subalgebra ${\cal C}_{\text{sp}}$,
as the representations do in the case of the translation 
action.  

 For any $u\in\text{Rep}(G)$, define the map
$$c_u: H_u\to\overline{L_u}\otimes{\cal C}$$ 
associated with
the spectral space $\overline{L_u}$ by:
$$c_u(\psi):=\sum_k\overline{T_k}\otimes T_k(\psi),\eqno(2.3)$$
where $T_k$ is any orthonormal basis of $L_u$.
Clearly $c_u$  is determined by its coefficients
$$c^u_{\overline{T},\psi}:=\overline{T}^*\otimes I 
c_u(\psi)=T(\psi),\quad\psi\in 
H_u,\  T\in L_u.\eqno(2.4)$$
The  $c_u$'s are  called {\it multiplicity 
maps} in \cite{PR}. 

In the example of the translation action, identifying
$\overline{L_u}$ and $H_u$ we have 
$c_u=u$.

We can represent $c_u$ as a rectangular matrix whose $k$-th
row is given by the multiplet $T_k=(T_k(\psi_1)\dots T_k(\psi_d))$
transforming like $u$ under $\delta$.

It is known that the set of all coefficients 
 $\{c^u_{\overline{T_k},\psi_j}=T_k(\psi_j), j,k\}$, 
of the 
multiplicity maps in orthonormal bases forms
a linear basis for the dense $^*$--subalgebra 
${\cal C}_{\text{sp}}$ \cite{BRV}, \cite{PR}, 
when $u$ varies in a complete set of irreducible representations of 
$\text{sp}(\delta)$,  
generalizing a well known 
property of 
matrix coefficients of a compact quantum group 
\cite{WLesHouches}.

In \cite{BRV} the authors introduce a new numerical invariant, the {\it 
quantum 
multiplicity} $m(u)$ of the representation $u$, in the following way.
If $j:H_u\to H_{\overline{u}}$ defines a conjugate representation of $u$
in the sense recalled in the previous subsection, then we can associate 
 an invertible 
antilinear 
$J:L_u\to L_{\overline{u}}$ with $J$ by setting $J(T)(\phi):=T(j^{-1}(\phi))^*$.
Its inverse $J^{-1}:L_{\overline{u}}\to L_u$ is given
by $J^{-1}(S)(\psi)=S(j(\psi))^*$. If $u$ is irreducible,
$m(u)^2:=\text{Trace}(JJ^*)\text{Trace}((JJ^*)^{-1}).$
One has:
$$\text{mult}(u)\leq m(u)\leq d(u),$$
an inequality which strengthens the inequality $\text{mult}(u)\leq d(u)$ 
previously obtained by Boca \cite{Boca}, when generalizing the HLS theorem 
\cite{HLS}
$\text{mult}(u)\leq \text{dim}(u)$ in the group case.
If $u$ is reducible, we define $m(u)$ as  the infimum of all 
the above trace values, derived from all possible solutions of the 
conjugate equations for $u$. Then the  inequality 
$$\text{dim}(L_u)\leq m(u)\leq 
d(u)\eqno(2.5)$$
holds for all representations $u$.
Notice that $m(u)$ takes the smallest possible value $\text{dim}(L_u)$ 
precisely when for some $j$ the associated $J$ is a scalar multiple of an 
antiunitary. 
Examples of ergodic actions of $S_\mu U(2)$ where 
$\text{dim}(u)<\text{mult}(u)<m(u)=d(u)$ have been constructed in 
\cite{BRV}.
\medskip

\noindent{\bf 2.3 The spectral functor of an ergodic action and quasitensor functors}
\medskip

It has been shown in \cite{PR} that ergodic actions of  compact 
quantum 
groups on unital $C^*$--algebras have a duality theory resembling the 
duality theory of Woronowicz for compact quantum groups: 
an ergodic $G$--action on ${\cal C}$ has a dual object allowing one to  
reconstruct the $G$--action on the maximal completion of 
 ${\cal C}_{\text{sp}}$. Furthermore, 
the dual objects of ergodic actions have been characterized. 

 The map $u\mapsto\overline{L_u}$ 
can be extended to a {\it functor} 
$$\overline{L}:\text{Rep}(G)\to{\cal H}$$
from the category  of representations of $G$ 
to the category ${\cal H}$ of Hilbert 
spaces. This functor is defined on arrows as follows.

If $A\in(u,v)$ and $T\in L_v$ then $T\circ A: H_u\to{\cal C}$
lies in $L_u$.
Hence 
if we identify 
$\overline{L_u}$ canonically with the dual vector space of $L_u$, any 
arrow
$A\in(u,v)$ in $\text{Rep}(G)$ induces a linear map
$\overline{L}_A\in(\overline{L_u},\overline{L_v})$
via the natural pairing between $L_u$ and $\overline{L_u}$: 
$$\overline{L}_A:\varphi\in\overline{L}_u\to (T\in L_v\to \varphi(T\circ 
A))\in\overline{L_v}.$$

The spectral functor $\overline{L}$ and the multiplicity maps $c_u$ 
are related as follows
$$\overline{L}_A\otimes I\circ c_u=c_v\circ A,\quad A\in(u,v),$$
for any $u,v\in\text{Rep}(G)$.
In terms of the matrix coefficients of $c_u$
this  reads
$$c^u_{\overline{L}_{A^*}\,\overline{S},\psi}=c^v_{\overline{S},A\psi},\quad 
A\in(u,v), \psi\in H_u, S\in L_v.\eqno(2.6)$$

Taking  the tensor $C^*$--category  structure of 
$\text{Rep}(G)$ and ${\cal H}$ into account  one can see that $\overline{L}$ is a 
$^*$--functor, but {\it not} a tensor $^*$--functor, in general.

 In fact, for 
$u,v\in\text{Rep}(G)$, the 
tensor product Hilbert 
space 
$\overline{L_u}\otimes\overline{L_v}$ 
is in general just a subspace
of $\overline{L_{u\otimes v}}$, in the sense that there is a 
natural isometric inclusion
$$\tilde{\overline L}_{u,v}:\overline{L_u}\otimes\overline{L_v}\to\overline{L_{u\otimes 
v}}$$
identifying a simple tensor $\overline{S}\otimes\overline{T}$ with
the complex conjugate of the element of $L_{u\otimes v}$ defined by
$$\psi\otimes\phi\in H_u\otimes H_v\to S(\psi)T(\phi).$$
The main result of \cite{PR} characterizes the set of all 
ergodic 
action duals 
$(\overline{L},\tilde{\overline L})$ algebraically among all
$^*$--functors
$${\tau}:\text{Rep}(G)\to{\cal H}$$
endowed  with isometries 
$\tilde\tau_{u,v}:{\tau}_u\otimes{\tau}_v\to{\tau}_{u\otimes v}$. 
These are precisely the {\it quasitensor functors}, defined below.

We shall refer to \cite{DRInventiones} for the notion of an abstract 
(strict) tensor 
$C^*$--category ${\cal T}$. The tensor product between objects 
$u$ and 
$v$ will be  denoted by $u\otimes v$ 
and between 
arrows $S$ and $T$ by $S\otimes T$. 
The tensor unit object will be denoted $\iota$.
We shall   assume that $\iota$ is an irreducible 
object: $(\iota,\iota)={\mathbb C}$, unless otherwise specified. 
The n-th tensor power of an object $u$ will be denoted $u^n$.
When we refer to a tensor $C^*$--category of Hilbert spaces, 
we mean that the objects are finite dimensional Hilbert spaces 
and contain Hilbert spaces of any finite dimension. The spaces 
of arrows are all linear operators between the Hilbert spaces in question.\medskip 

Let   ${\cal T}$ and ${\cal R}$ be strict tensor $C^*$--categories.
A $^*$--functor ${\tau}: {\cal
T}\to {\cal R}$ together with a collection of isometries
$\tilde\tau_{u,v}\in(\tau_u\otimes\tau_v, \tau_{u\otimes v})$,
 for objects $u$, $v\in{\cal T}$,
is called
{\it quasitensor} if
$${\tau}_\iota=\iota,\eqno(2.7)$$
$$\tilde\tau_{u,\iota}=\tilde\tau_{\iota, u}=1_{\tau_u},\eqno(2.8)$$
$$\tilde\tau_{u,v\otimes w}^*\circ\tilde\tau_{u\otimes v,w}=1_{\tau_u}\otimes\tilde\tau_{v,w}
\circ{\tilde\tau_{u,v}}^*\otimes 1_{\tau_w}\eqno(2.9)$$ 
and if
$$\tau({S\otimes T})\circ
\tilde\tau_{u,v}=\tilde\tau_{u',v'}\circ\tau(S)\otimes\tau(T),\eqno(2.10)$$
 for any other pair of objects $u'$, $v'$ and arrows
$S\in(u, u')$,
$T\in(v,v')$. 
In particular, a 
{\it tensor} functor $\tau$ is quasitensor with $\tilde\tau_{u,v}:=1_{\tau_u\otimes\tau_v}$
, as
$\tau_{u\otimes v}=\tau_{u}\otimes\tau_{v}$ for all objects $u$, $v$,   
$(2.7)$ and $(2.8)$ hold by assumption and
$(2.9)$ and $(2.10)$ are trivially satisfied.  
More  generally, if all the isometries $\tilde\tau_{u,v}$ are unitary, 
we recover the known notion of a {\it relaxed tensor} functor.

This definition of a quasitensor functor differs from that given in 
\cite{PR} and the equivalence is established in the appendix.\medskip 

Given a quasitensor functor $(\mu,\tilde\mu)$ into the category of Hilbert  
spaces, let $\tau$ denote the embedding functor of the category of finite dimensional unitary representations of a compact quantum group $G$ into the category of Hilbert spaces then, as shown in \cite{PR}, there is a canonical ergodic action of $G$ on a $C^*$--algebra ${}_\mu\cC_\tau$. 
If $\mu$ is the spectral functor of an ergodic action of $G$ on a $C^*$--algebra $\cB$ then $\mu$ is isomorphic to the spectral functor of the 
action on ${}_\mu\cC_\tau$ and the dense spectral subalgebras of 
$\cB$ and ${}_\mu\cC_\tau$ are canonically isomorphic. However, $\cB$ and ${}_\mu\cC_\tau$ need not be isomorphic. ${}_\mu\cC_\tau$ is the 
completion of its dense spectral subalgebra in the maximal $C^*$--norm 
and this may not be the case for $\cB$.\medskip

\medskip

\end{section}

\begin{section} {The main results}

In this section we state our main results. \medskip

\noindent{\bf 3.1 Theorem} {\sl Let $x$ be an object of a tensor 
$C^*$--category with irreducible tensor unit $\iota$ and let $R\in(\iota,x^2)$ satisfy $R^*\otimes 1_x\circ 1_x\otimes R=\pm 1_x$
and $\|R\|^2\geq2$.
For any integer $2\leq n\leq\|R\|^2$, let
$Q\in M_n({\mathbb C})$ be any invertible matrix satisfying
$$Q\overline{Q}=\pm I,\quad \text{Trace}(Q^*Q)=\text{Trace}(Q^*Q)^{-1}=\|R\|^2.$$
Then there is an ergodic action of the compact quantum group $B_u(Q)$ of Wang \cite{Wang2} on a unital
$C^*$--algebra ${\cal C}$ with spectral spaces
$\overline{L}_{u^r}=(\iota,x^r)$, 
$r\geq0$ and
$\overline L_{\sum_k\psi_k\otimes 
Q^*\psi_k}= R$, where the sum is taken over an orthonormal basis $\psi_k$ of the Hilbert 
space of $u$.
 If $m(u)={\rm dim}(\iota,x)$, $u$ being the defining representation of $B_u(Q)$, then 
$$m(u^r)=\text{dim}(\iota,x^r).$$ 
In particular, if we choose $n=2$ we get an ergodic action of $S_{\mu} U(2)$ for a nonzero
$-1<\mu<1$ determined by $|\mu+\mu^{-1}|=\|R\|^2$, where $\mu>0$ if and only if $x$ is pseudoreal.}\medskip 

 In the examples derived from subfactors and treated in \cite{PR1}. We do 
 have $m(u)={\rm dim}(\iota, x)$.\medskip

\noindent
{\bf 3.2 Theorem} {\sl Let $x$ be an object of a tensor 
$C^*$--category with irreducible tensor unit $\iota$ and let $R\in(\iota,\overline x\otimes x)$ and $\overline R\in(\iota,x\otimes\overline x)$ satisfy $R^*\otimes 1_{\overline x}\circ 1_{\overline x}\otimes R=1_{\overline x}$ and $\overline R^*\otimes 1_x\circ 1_x\otimes R=1_x$
and $\|R\|^2=\|\overline R\|^2\geq2$. For any integer $2\leq n\leq\|R\|^2$, let
$Q\in M_n({\mathbb C})$ be any positive invertible matrix satisfying
$$\text{Trace}(Q)=\text{Trace}(Q^{-1})=\|R\|^2.$$
Then there is an ergodic action of the compact quantum group $A_u(Q)$ of Wang on a unital $C^*$--algebra ${\cal C}$ with spectral spaces
$ \overline L_{q(u,\overline u)}=(\iota,q(x,\overline x))$, 
where $q$ is a monomial in two variables and $u$ the defining representation of $A_u(Q)$. If $m(u)={\rm dim}(\iota,x)$, 
$$m(q(u,\overline u))={\rm dim}(\iota,q(x,\overline x)),$$ for each $q$.
}\medskip

The proofs involve two main steps. The first is to embed the tensor $^*$-- subcategory generated by $R$ or by $R$ and $\overline R$ into the category of Hilbert spaces. The second step is 
to define the ergodic  action, by applying the duality theorem for ergodic actions of compact quantum groups on unital $C^*$--algebras proved in \cite{PR}. 
The construction of the $C^*$--algebra will be given 
in Sect. 4 in greater generality than in \cite{PR} and the $G$--action is explained in Sect. 5. We start 
with {\it a 
pair} of quasitensor functors $(\tau,\tilde{\tau})$, $(\mu,\tilde{\mu})$ and obtain the following result.
\medskip

\noindent{\bf 3.3 Theorem} {\sl
Let $\cA$ be  a tensor $C^*$--category  with conjugates and $(\mu,\tilde\mu):\cA\to{\cal M}$ 
and $(\tau,\tilde\tau):\cA\to\cT$ quasitensor functors. We let $^\circ_\mu\cC_\tau$ be the linear space 
$\sum_u(\mu_u,\iota)\otimes(\iota,\tau_u)$, the sum being taken over the objects of $\cA$, 
quotiented by the linear subspace generated by elements of the form 
$$M\circ\mu(A)\otimes T-M\otimes\tau(A)\circ T.$$ 
Thus we may write $^\circ_\mu\cC_\tau=\sum_u(\mu_u,\iota)\otimes_\cA(\iota,\tau_u)$.   
Then 
\begin{description}
\item{a)}
$^\circ_\mu\cC_\tau$ can be given the structure of a $^*$--algebra. 
\item{b)}
If either $\mu$ or $\tau$ is tensorial then    $^\circ_\mu\cC_\tau$ can be completed in the maximal $C^*$--norm to give a unital $C^*$--algebra ${}_\mu\cC_\tau$. 
\item{c)}
If we have an action $\eta$ of a compact qantum group $G$ on ${\cal T}$ leaving the objects invariant, there is 
a unique action $\alpha$ of $G$ on ${}_\mu\cC_\tau$ such that $\alpha(M\otimes T)=M\otimes\eta(T)$.
\end{description}

}\medskip

\end{section}

\begin{section}{$C^*$--algebras from pairs of quasitensor functors}

Let $\cA$ be  a tensor $C^*$--category with conjugates and $(\mu,\tilde\mu):\cA\to{\cal M}$ 
and $(\tau,\tilde\tau):\cA\to\cT$ quasitensor functors. We let $^\circ_\mu\cC_\tau$ be the linear space 
$\sum_u(\mu_u,\iota)\otimes(\iota,\tau_u)$, the sum being taken over the objects of $\cA$, 
quotiented by the linear subspace generated by elements of the form 
$$M\circ\mu(A)\otimes T-M\otimes\tau(A)\circ T.$$ 
Thus we may write $^\circ_\mu\cC_\tau=\sum_u(\mu_u,\iota)\otimes_\cA(\iota,\tau_u)$. We next 
define a product on $^\circ_\mu\cC_\tau$ setting for $L\in(\mu_u,\iota)$, $M\in(\mu_v,\iota)$, $S\in(\iota,\tau_u)$,
$T\in(\iota,\tau_v)$,
$$(L\otimes S)(M\otimes T):=(L\otimes M)\circ\tilde\mu_{u,v}^*\otimes\tilde\tau_{u,v}\circ(S\otimes T).$$
It is easy to check that the product is well defined and associative. 

When either $(\tau,\tilde\tau)$ or $(\mu,\tilde\mu)$ is minimal  in the sense defined in the appendix $^\circ_\mu\cC_\tau$ reduces to the complex numbers. The reason is that this algebra does not change if we complete $\cA$ under 
direct sums and subobjects and extend $(\tau,\tilde\tau)$. Every object of $\cA$ is then a direct sum 
of irreducibles and it becomes clear that we can restrict the sum over $u$ in the definition of $^\circ_\mu\cC_\tau$ to a representative set of irreducibles. But then $(\iota,\mu_u)=0$ unless $u=\iota$ so $^\circ_\mu\cC_\tau=(\iota,\iota).$

Tensor $C^*$--categories with 
conjugates have been studied in \cite{LR}. We recall the notion of 
conjugate object $\overline{u}$ of $u$. This object is defined, up to 
unitary equivalence, by 
the existence of two intertwiners $R\in(\iota,\overline{u}\otimes u)$,
$\overline{R}\in(\iota, u\otimes\overline{u})$ satisfying the conjugate 
equations:
$$\overline{R}^*\otimes 1_u\circ 1_u\otimes R=1_u,\eqno(4.1)$$
$$R^*\otimes 1_{\overline{u}}\circ 
1_{\overline{u}}\otimes\overline{R}=1_{\overline{u}}.\eqno(4.2)$$
The {\it intrinsic dimension} $d(u)$ of $u$ is 
the infimum of all possible $\|R\|\|\overline{R}\|$.

If $G$ is a compact quantum group,  $\text{Rep}(G)$   is a tensor $C^*$--category with conjugates: for any representation $u$ with conjugate representation $\overline{u}$ defined by the 
antilinear intertwiner $j:H_u\to H_{\overline{u}}$ as in $(2.2)$, the elements
$R:=\sum\psi_j\otimes j^{-1}\psi_j$ and
$\overline{R}:=\sum_k\phi_k\otimes j\phi_k$ are intertwiners in 
$(\iota, \overline{u}\otimes 
u)$ and  $(\iota, 
u\otimes\overline{u})$ respectively
and satisfy the conjugate equations.
Hence every representation has an associated  intrinsic dimension 
$d(u)$, also called quantum dimension,
given by 
$$d(u)^2=\text{inf}(\|R\|\|\overline{R}\|)^2=
\text{inf}\,\text{Trace}(j^*j)\text{Trace}((j^*j)^{-1}).\eqno(4.3)$$
Notice that $d(u)\geq\text{dim}(u)$ with equality 
if and only if $j$ is antiunitary.
In terms of the quantum group, the condition $\text{dim}(u)=d(u)$ for 
all $u$ is equivalent to  
requiring the coinverse $\kappa$ to be involutive.

Let $\cA$ be a tensor $C^*$--category and pick for each object $u$ of $\cA$ a solution $R_u,\overline R_u$ of the conjugate equations. We agree to take $R_\iota=\overline R_\iota=1_\iota$. 
There is an associated conjugation on $\cA$ defined, for $A\in(v,u)$, by 
$$A^\bullet:=R^*_v\otimes 1_{\overline u}\circ 1_{\overline v}\otimes A^*\otimes 1_{\overline u}\circ 
1_{\overline v}\otimes\overline R_u.$$ 
$A^\bullet\in(\overline v,\overline u)$ can also be defined by the equation 
$$1_v\otimes A^\bullet\circ \overline R_v=A^*\otimes 1_{\overline u}\circ\overline R_u.$$
If $B\in(w,v)$ then $(A\circ B)^\bullet=A^\bullet\circ B^\bullet$. If we use the product solutions of the 
conjugate equations for defining the conjugate of a product then $(A\otimes B)^\bullet=B^\bullet\otimes A^\bullet$. In fact, if $A\in(u,u')$ and $B\in(v,v')$ then
$$(A\otimes B)^\bullet=R^*_{u\otimes v}\otimes 1_{\overline v'\otimes\overline u'}\circ 1_{\overline v\otimes\overline u}\otimes A^*\otimes B^*\otimes 1_{\overline v'\otimes\overline u'}\circ 1_{\overline v\otimes\overline u}\otimes \overline R_{u'\otimes v'}.$$ 
Substituting in the product form of the solutions we get 
$$(A\otimes B)^\bullet=$$ 
$$(R^*_v\circ 1_{\overline v}\otimes R^*_u\otimes 1_v)\otimes 1_{\overline v'\otimes\overline u'}\circ 1_{\overline v\otimes\overline u}\otimes A^*\otimes B^*\otimes 1_{\overline v'\otimes\overline u'}\circ 1_{\overline v\otimes\overline u}\otimes(1_{u'}\otimes\overline R_{v'}\otimes 1_{\overline u'}\circ\overline R_{u'})=$$ 
$$(R^*_v\circ 1_{\overline v}\otimes B^*\circ 1_{\overline v}\otimes R^*_u\otimes 1_{v'})\otimes 1_{\overline v'\otimes\overline u'}\circ 1_{\overline v\otimes\overline u}\otimes(1_u\otimes\overline R_{v'}\otimes 1_{\overline u'}\circ  A^*\otimes 1_{\overline u'}\circ\overline R_{u'})=$$ 
 $$(R^*_v\circ 1_{\overline v}\otimes B^*)\otimes 1_{\overline v'\otimes\overline u'}\circ 1_{\overline v}\otimes\overline R_{v'}\otimes 1_{\overline u'}\circ 1_{\overline v}\otimes R^*_u\otimes 1_{\overline u'}\circ 1_{\overline v\otimes\overline u}\otimes(A^*\otimes 1_{\overline u'}\circ\overline R_{u'})=$$
 $$B^\bullet\otimes 1_{\overline u'}\circ 1_{\overline v}\otimes A^\bullet=B^\bullet\otimes A^\bullet.$$
 A computation shows that the inverse of $A\to A^\bullet$ is 
$$A=\overline R_v^*\otimes 1_u\circ 1_v\otimes A^{\bullet*}\otimes 1_u\circ 1_v\otimes R_u.$$
Now 
$$\hat R_u:=\tilde\mu_{\overline u,u}^*\circ\mu(R_u),\hat{\overline R_u}:=\tilde\mu_{u,\overline u}^*\circ\mu(\overline R_u)\eqno(4.4)$$ 
is a solution of the conjugate equations for $\mu_u$ since 
$$\hat R_u^*\otimes 1_{\mu_{\overline u}}\circ 1_{\mu_{\overline u}}\otimes \hat{\overline R}_u=\mu(R_u^*)\otimes 1_{\mu_{\overline u}}\circ\tilde\mu_{\overline u,u}\otimes 1_{\mu_{\overline u}}\circ 1_{\mu_{\overline u}}\otimes\tilde\mu_{u,\overline u}^*\circ 1_{\mu_{\overline u}}\otimes\mu(\overline R_u)=$$
$$\mu(R_u^*)\otimes 1_{\mu_{\overline u}}\circ\tilde\mu_{\overline u\otimes u,\overline u}^*\circ\tilde\mu_{\overline u,u\otimes\overline u}\circ 1_{\mu_{\overline u}}\otimes\mu(\overline R_u)=\mu(R_u^*\otimes 1_{\overline u})\circ\mu(1_{\overline u}\otimes\overline R_u)=1_{\mu_{\overline u}}$$ 
with the other relation following similarly. There is therefore a conjugation defined on the full subcategory of $\cM$ whose objects are the images of those of $\cA$ under $\mu$.
\medskip

\noindent
{\it Remark} If $\mu$ is not injective on objects, this conjugation is not well defined. This plays no role in the following since $\hat R_u$ is labelled by $u$ rather than $\mu_u$. In what follows, the tensor category $\cM$ can, if desired, be replaced by a tensor $C^*$--category whose objects are those of $\cA$ and where the arrows from $u$ to $v$ are arrows from $\mu_u$ to $\mu_v$ with the obvious algebraic operations. The $^*$--functor $\mu$ then becomes an isomorphism on objects.\medskip 
 
Changing the solution of the conjugate equations using an invertible $X$ (see Appendix), then as we would expect, the corresponding change in $\hat R_u$ and $\hat{\overline R}_u$ is induced by 
$\mu(X)$. The solutions of the conjugate equations for $\tau_u$, defined analogously, will be 
denoted by $\tilde R_u,\tilde{\overline R}_u$.

Given $A\in(v,u)$, then 
$$\mu(A^\bullet)=\mu(R^*_v\otimes 1_{\overline u}\circ 1_{\overline v}\otimes A^*\otimes 1_{\overline u}\circ 1_{\overline v}\otimes\overline R_u)=$$
$$\mu(R^*_v)\otimes 1_{\mu_{\overline u}}\circ\tilde\mu_{\overline v\otimes v,\overline u}^*\circ\tilde\mu_{\overline v,v\otimes\overline u}\circ 1_{\mu_{\overline v}}\otimes\mu(A^*\otimes 1_{\overline u}\circ\overline R_u)=$$
$$\mu(R_v^*)\otimes 1_{\mu_{\overline u}}\circ\tilde\mu_{\overline v,v}\otimes 1_{\mu_{\overline v,v}}\circ\tilde1_{\mu_{\overline v}}\circ\tilde\mu_{v,\overline u}^*\circ 1_{\mu_{\overline v}}\otimes\mu(A^*\otimes 1_{\overline u}\circ \overline R_u)=$$
$$\hat R_v^*\otimes 1_{\mu_{\overline u}}\circ 1_{\mu_{\overline v}}\otimes\mu(A^*)\otimes 1_{\mu_{\overline u}}\circ 1_{\mu_{\overline v}}\otimes\tilde\mu_{u,\overline u}^*\circ 1_{\mu_{\overline v}}\otimes \overline R_u=$$
$$\hat R_v^*\otimes 1_{\mu_{\overline u}}\circ  1_{\mu_{\overline v}}\otimes\mu(A^*)\otimes 1_{\mu_{\overline u}}\circ 1_{\mu_{\overline u}}\otimes\hat {\overline R}_v. $$ 
Thus $\mu(A^\bullet)=\mu(A)^\bullet$.\medskip 

\noindent

When $(\mu,\tilde\mu)$ and $(\tau,\tilde\tau)$ are quasitensor $^*$--functors, we define an involution on $^\circ_\mu\cC_\tau$ by setting 
$$(M\otimes T)^*:=M^\bullet\otimes T^\bullet.$$
This is well defined since, for example, 
$$(M\circ\mu(A)\otimes T)^*=(M\circ \mu(A))^\bullet\otimes T^\bullet=M^\bullet\circ\mu(A^\bullet)\otimes T^\bullet=$$
$$M^\bullet\otimes\tau(A^\bullet)\circ T^\bullet=(M\otimes\tau(A)\circ T)^*.$$
If we change the solution of the conjugate equations using an invertible $X$ (see Appendix), then $(M\otimes T)^*$ 
becomes $M^\bullet\circ \mu(X^*)\otimes\tau(X^{-1*})\circ T^\bullet=(M\otimes T)^*$. In other words, 
the involution is independent of the choice of solutions of the conjugate equations in $\cA$. Thus to 
check that we really have an involution, it suffices to pick $R_{\overline u}=\overline R_u$ and 
$\overline R_{\overline u}=R_u$ when computing the second adjoint. In this case, the above 
computation of the inverse of $A\mapsto A^\bullet$ implies that we have an involution.\medskip 

\noindent
{\bf 4.1  Proposition} {\sl The product and involution defined above make 
$^\circ_\mu\cC_\tau$ into a $^*$--algebra.}\medskip 

\noindent
{\it Proof} It suffices to show that 
$$(N^\bullet\otimes M^\bullet\circ\tilde\mu_{\overline v,\overline u}^*)\otimes(\tilde\tau_{\overline v,\overline u}\circ T^\bullet\otimes S^\bullet)=(N^\bullet\otimes T^\bullet)\otimes(M^\bullet\otimes S^\bullet)=$$
$$(M\otimes N\circ\tilde\mu_{u,v}^*)^\bullet\otimes(\tilde\tau_{u,v}\circ S\otimes T)^\bullet.$$ 
As the involution 
is independent of the choice of solutions of the conjugate equations, we may suppose, in evaluating 
this expression, that $\hat R_{u\otimes v}=\tilde\mu_{\overline v\otimes\overline u,u\otimes v}^*\circ\mu(1_{\overline v}\otimes R_u\otimes 1_v\circ R_v)$ with an analogous expression for $\tilde R_{u\otimes v}$. Now 
$$(M\otimes N\circ\tilde\mu_{u,v}^*)^\bullet=\hat R_{u\otimes v}^*\circ 1_{\mu_{\overline v\otimes\overline  
u}}\otimes\tilde\mu_{u,v}\circ 1_{\mu_{\overline v\otimes\overline  u}}\otimes M^*\otimes N^*= $$ 
$$\mu(R_v^*)\circ\mu(1_{\overline v}\otimes R_u^*\otimes 1_v)\circ\tilde\mu_{\overline v\otimes\overline u,u\otimes v}\circ 1_{\mu_{\overline v\otimes\overline u}}\otimes\tilde\mu_{u,v}\circ 1_{\mu_{\overline v\otimes\overline  u}}\otimes M^*\otimes N^*= $$ 
$$\mu(R_v^*)\circ\mu(1_{\overline v}\otimes R_u^*\otimes 1_v)\circ\tilde\mu_{\overline v\otimes\overline
u\otimes u, v}\circ\tilde\mu_{\overline v\otimes\overline u,u}\otimes 1_{\mu_v}\circ 1_{\mu_{\overline v\otimes\overline  u}}\otimes M^*\otimes N^*= $$ 
$$\mu(R_v^*)\circ\tilde\mu_{\overline v,v}\circ 1_{\mu_{\overline v}}\otimes\mu(R^*_u)\otimes 1_{\mu_v}
\circ\tilde\mu_{\overline v,\overline u\otimes u}^*\otimes 1_{\mu_v}\circ\tilde\mu_{\overline v\otimes\overline u,u}\otimes 1_{\mu_v}\circ 1_{\mu_{\overline v\otimes\overline  u}}\otimes M^*\otimes N^*= $$ 
$$\hat R_v^*\circ 1_{\mu_{\overline v}}\otimes N^*\circ 1_{\mu_{\overline v}}\otimes\mu(R^*_u)\circ
1_{\mu_{\overline v}}\otimes\tilde\mu_{\overline u,u}\circ\tilde\mu_{\overline v,\overline u}\otimes 1_{\mu_u}\circ 1_{\mu_{\overline v\otimes\overline u}}\otimes M^*=$$ 
$$(\hat R^*_v\circ 1_{\mu_{\overline v}}\otimes N^*)\otimes(\hat R_u^*\circ 1_{\mu_{\overline u}}\otimes M^*)\circ\tilde\mu_{\overline v,\overline u}^*=(N^\bullet\otimes M^\bullet\circ\tilde\mu_{\overline v,\overline u}^*).$$ 
This proves the result since the term involving $S$ and $T$ can be treated in the same way.\medskip 

\noindent
{\it Remark} Note that $^\circ_\mu\cC_\tau$ depends only on the images of $\mu$ and $\tau$. However 
the images will not in general be tensor categories and the existence of $^\circ_\mu\cC_\tau$ depends 
on having two quasitensor functors.\medskip 

In the following $c$ will denote the support of $\iota$ (see Appendix).\medskip

\noindent
{\bf 4.2 Corollary } {\sl If $\tilde\tau_{u,v}^\bullet\circ c_{\tau_{\overline v}}\otimes c_{\tau_{\overline u}}=\tilde\tau_{u,v}^{*\bullet*}\circ c_{\tau_{\overline v}}\otimes c_{\tau_{\overline u}}=\tilde\tau_{\overline v,\overline u}\circ c_{\tau_{\overline v}}\otimes c_{\tau_{\overline u}}$.}\medskip 

\noindent
{\it Proof.} Let $S_i\in(\iota,\tau_u)$, $T_j\in(\iota,\tau_v)$ be orthonormal bases then in the proof  of Proposition 4.1 we have seen that 
$$\tilde\tau_{u,v}^\bullet\circ T_j^\bullet\otimes S_i^\bullet=(\tilde\tau_{u,v}\circ S_i\otimes T_j)^\bullet=
\tilde\tau_{\overline v,\overline u}\circ T_j^\bullet\otimes S_i^\bullet.$$
Multiplying on the right by $T_j^{*\bullet}\otimes S_i^{*\bullet}$, summing over $i$ and $j$ and using the fact that, by Lemma A.5, $c_{\tau_u}^\bullet=c_{\tau_{\overline u}}$, we get the one equality. The other equality follows similarly from Proposition 4.1 but using the part involving the functor $\mu$.\medskip

We now make use of the irreducibility of the tensor unit $\iota$  in $\cA$, 
$\cM$ and $\cT$ and define a linear functional $h$ on $^\circ_\mu\cC_\tau$ by setting for $M\in(\mu_u,\iota)$ and $T\in(\iota,\tau_u)$, 
$$h(M\otimes T):=(M\circ\mu(c_u))\otimes T=\sum_i(M\circ\mu(V_{u,i}))\otimes(\tau(V_{u,i}^*)\circ T),$$ 
where $V_{u,i}$ is an orthonormal basis in $(\iota,u)$. 
The first expression shows that $h$ is well defined and the second that it takes values in $\mathbb C$. 
\medskip

The next task is to show that $h$ is a faithful positive linear functional and it would be natural to argue in terms of irreducibles. However $\cA$ does not necessary have sufficient irreducibles and there are two ways 
to proceed. Let $\cB$ denote the completion of $\cA$ under subobjects then $\cB$ has sufficient 
irreducibles and we have a canonical inclusion functor from $\cA$ to $\cB$. The quasitensor functors 
$(\mu,\tilde\mu)$ and $(\tau,\tilde\tau)$ from $\cA$ can be extended to quasitensor functors $(\nu,\tilde\nu)$ and $(\sigma,\tilde\sigma)$ from $\cB$. A variant of Proposition 4.4 below shows that ${}^\circ_\mu\cC_\tau$ and ${}^\circ_\nu\cC_\sigma$ are canonically isomorphic. After this we may suppose that $\cA$ has sufficient irreducibles. To avoid giving the details involved, we give an alternative proof using 
minimal projections in $\cA$ in place of irreducibles, every unit being a sum of minimal projections.\medskip 

As we have seen that the involution on our algebra is independent of the choice of conjugate, we suppose 
in the following computation that $u\mapsto R_u$ is a standard choice of solutions of the conjugate equations (see Appendix). The conjugation then commutes with the adjoint and maps projections into 
projections. 
If $E\in(u,u)$ is a minimal projection then $E^\bullet$ is a minimal projection in $(\overline u,\overline u)$. Setting $R_E:=E^\bullet\otimes E\circ R_u$ and $\overline R_E:=E\otimes E^\bullet\circ\overline R_u$, 
we have 
$$E\otimes R^*_E\circ\overline R_E\otimes E=E,\quad E^\bullet\otimes\overline R^*_E\circ R_E\otimes E^\bullet=E^\bullet,$$ 
the form taken by the conjugate equations for minimal projections. If $E$ is any projection in $(u,u)$ we 
let $V_{E,i}$, $i=1,2,\dots n_E$, be a maximal set of mutually orthogonal isometries in $(\iota,u)$ with $V_{E,i}\circ V_{E.i}^*\leq E$ and we set $c_E:=\sum_iV_{E,i}\circ V_{E,i}^*$. 
Two minimal projections are equivalent if they are connected by a partial isometry and we pick a set $\hat E$ of minimal projections, one from each equivalence class.
Note that $c_{E^\bullet\otimes F}=0$ if $E$ and $F$ are inequivalent minimal projections whereas 
$c_{E^\bullet\otimes E}=\|R_E\|^{-2}R_E\circ R_E^*$.\medskip

\noindent
{\bf 4.3 Proposition} {\sl $h$ is a faithful positive linear  functional 
on $^\circ_\mu\cC_\tau$.}\medskip 

\noindent
{\it Proof} Given an object $u$ of $\cA$ there are partial isometries $U_i$ with $U_i^*\circ U_i\in\hat E$ and $\sum_iU_i\circ U_i^*=1_u$. If $M\in(\mu_u,\iota)$ and $T\in(\iota,\tau_u)$ then 
$M\otimes T=\sum_iM\circ\mu(U_i)\otimes\tau(U_i^*)\circ T).$
Thus any element of $^\circ_\mu\cC_\tau$ is a sum of elements of the form 
$M\otimes T$ where $M=M\circ\mu(E)$ and $\tau(E)\circ T$ for some $E\in\hat E$. Given  
$L\otimes S$ with $L=L\circ\mu(F)$ and $S=\tau(F)\circ S$ then by the remarks above, 
$h((L\otimes S)^*(M\otimes T))=0$ if $E\neq F$ whereas if $E=F$, 
$$h((L\otimes S)^*(M\otimes T))=\|R_E\|^{-2}(\hat R_E^*\circ \mu(E^\bullet)\otimes L^*\otimes M\circ\hat R_E)(\hat R_E^*\circ\mu(E^\bullet)\otimes S^*\otimes T\circ\hat R_E)$$ 
$$=\|R_E\|^{-2}(\phi_E(L^*\circ M))(\phi_E^*(S^*\circ T)),$$ 
with $\phi_E$ the scalar product on $\{X:X\circ\mu(E)=X=\mu(E)\circ X\}$ and $\{Y:Y\circ\tau(E)=Y=\tau(E)\circ Y\}$ induced by $\hat{R_E}:=\mu(E^\bullet)\otimes\mu(E)\circ\hat{R}_u$ and $\tilde{R_E}:=\tau(E^\bullet)\otimes\tau(E)\circ\tilde R_u$ as in the appendix.
To complete the proof it is enough to show that $h(X^*X)\geq 0$ for $X:=\sum_{m,n}\lambda_{m,n}L_m\otimes S_n$ when $L_m$ is an orthonormal basis in $(\mu_u,\iota)\circ\mu(E)$, $S_n$ an orthonormal basis in $\tau(E)\circ(\iota,\tau_u)$ with respect to the scalar products $\phi_E$ and $\lambda_{m,n}\in\mathbb C$ 
with equality if and only if $\sum_m(L_m\otimes S_m)=0$. But the above computation shows that 
$$h(\sum_m(L_m\otimes S_m)^*\sum_n(L_n\otimes S_n))=$$
$$\sum_{m,n,p,q}\|R_E\|^{-2}\overline\lambda_{m,n}\lambda_{p,q}\phi_E(L_m^*\circ L_p)\phi_E(S_m^*\circ S_q)=\sum_{m,n}|\lambda_{m,n}|^2,$$ 
as required.\medskip 

The previous proposition implies in particular that $^\circ_\mu\cC_\tau$ has a non trivial $C^*$--norm. In general,   the maximal $C^*$--norm may not be finite. However, this will be the case  if either $\mu$ or $\tau$ is a (not necessarily strict) tensor  functor, and this suffices for the purposes of \cite{PRinduction}. To see this, we may   extend  the arguments    of \cite{PR}, where finiteness of the maximal $C^*$--norm    is explicitly shown,  to the setting of this paper.  However, the corresponding algebra ${\mathcal C}_{\mathcal F}$ of that paper  was introduced in a slightly different way,   by means of a complete set of irreducible objects. The following argument should make it clear that the two approaches are in fact equivalent.

 We  define a set of linear functionals on $^\circ_\mu\cC_\tau$. Pick a maximal set $E_k\in(u_k,u_k)$, $k\in K$, of inequivalent 
 minimal projections in $\cA$. Then for each $u$, pick partial isometries $W_i$, $i\in I_u$ such that 
 $W^*_i\circ W_i=E_{f_u(i)}$ where $f_u:I_u\to K$ and $\sum_{i\in I_u}W_i\circ W_i^*=1_u$. Given $k\in K$ and $M,M'\in(\iota,\mu_{u_k})$ and $N,N'\in(\iota,\mu_v)$, we set 
$$\omega_{T,M}(N^*\otimes S):=\sum_{i\in I_v,f_v(i)=k}(N^*\circ\mu(W_i)\circ M)(T^*\circ\tau(W_i^*)\circ S).$$ 
This expression is independent of the choice of the partial isometries and is understood to be zero if $f_v^{-1}(k)$ is the empty set. We must check that $\omega_{T,M}$ is well defined. To this end, let $A\in(v,w)$, and $P\in(\iota,\mu_w)$, then 
$$\omega_{T,M}(P^*\otimes\tau(A)\circ S)=\sum_{\ell\in I_w,f_w(\ell=k}(P^*\circ\mu(W_\ell)\circ M)(T^*\circ\tau(W_\ell^*)\circ\tau(A)\circ S)=$$
$$\sum_{i\in I_v,f_v=k,\ell\in K_v, f_w=k}(P^*\circ\mu(W_\ell)\circ M)(T^*\circ\tau(W_\ell^*)\circ\tau(A)\circ \tau(W_i)\circ\tau(W_i^*)\circ S)=$$
$$\sum_{i\in I_w;f_v (i)=k}(P^*\circ\mu(A)\circ\mu(W_i)\circ M)(T^*\circ\tau(W_i^*)\circ S)=\omega_{T,M}(P^*\circ\mu(A)\otimes S),$$ 
as required. Note that $\omega_{1_\iota,1_\iota}$ is just the Haar state $h$. One would expect 
$\omega_{M,M}$ to be a positive linear functional and imitating the proof in the case of $h$ should 
shed light on the question.\medskip

Since for $M,M'\in(\iota,\mu_u)$ and $W_i$, $i\in I_u$, as above, 
$$(M^*\otimes T)=\sum_{i\in I_u}(M^*\circ\mu(W_i))\otimes((\tau(W_i^*)\circ S),$$ 
every element of $^\circ_\mu\cC_\tau$ is a sum of elements of the form $(N^*\otimes S)$, 
where $\mu(E_k)\circ N=N$ and $\tau(E_k)\circ S=S$ for some $k\in K$. If $M=\mu(E_j)\circ M$ and $T=\tau(E_j)\circ T$, then 
$\omega_{T,M}(N^*\otimes S)=\delta_{jk}(T^*\circ S)(N^*\circ M)$. We now claim that the set of 
linear functionals $\omega_{T,M}$ with $M=\mu(E_k)\circ M$ and $T=\tau(E_k)\circ T$ separates sums of
elements of the form $N^*\otimes S$ where $N=\mu(E_k)\circ N$, $S=\tau(E_k)\circ S$. Any such sum $X$
may be written in the form $\sum_{i,j}\lambda_{ij}(M^*_i\otimes T_j)$, where $M_i$ and $T_j$ are orthonormal bases of the range of $\mu(E_k)$ and $\tau(E_k)$. $\omega_{M_n,M_p}(X)=\lambda_{np}$. Thus the above set of linear functionals  separates $^\circ_\mu\cC_\tau$.\medskip 

We let $_\mu\cC_\tau$ denote the $C^*$--algebra obtained by completing $^\circ_\mu\cC_\tau$ 
in the maximal $C^*$--norm.\medskip

We now investigate the functorial properties of the above construction. Let $(\eta,\tilde\eta):\cA_1\to\cA_2$, $(\mu_1,\tilde\mu_1):\cA_1\to{\cal M}$, $(\mu_2,\tilde\mu_2): \cA_2\to{\cal M}$, $(\tau_1,\tilde\tau_1):\cA_1\to\cT$ and $(\tau_2,\tilde\tau_2):\cA_2\to\cT$ be quasitensor functors with $(\mu_1,\tilde\mu_1)=(\mu_2,\tilde\mu_2)\circ(\eta,\tilde\eta)$ and $(\tau_1,\tilde\tau_1)=(\tau_2,\tilde\tau_2)\circ(\eta,\tilde\eta)$. 
Then the above equalities imply that there is a well defined natural unital multiplicative map $\eta_*$ from $_{\mu_1}\cC_{\tau_1}$ to $_{\mu_2}\cC_{\tau_2}$. Since the adjoint is independent of the choice of solutions of the conjugate equations we may suppose that if $R_u$ is chosen in $\cA_1$ then $\eta(R_u)$ is chosen in $\cA_2$ and a computation now shows that $\hat R_{\eta_u}=\hat R_u$ so that $\eta_*$ is a unital morphism. Obviously, $\eta\mapsto\eta_*$ is a covariant functor.\medskip

\noindent
{\bf  4.4 Proposition} {\sl If $(\eta,\tilde\eta):\cA_1\to\cA_2$ is full and each object of $\cA_2$ is a direct 
sum of the images of projections under $\eta$ then $\eta_*$ is an isomorphism.}\medskip

\noindent
{\it Proof} As the $C^*$--algebras in question are obtained by completion in the maximal $C^*$--norm 
it will suffice to show that we have an isomorphism before completion. 
Given $M\in(\iota,\mu_{2_,x})$ and $T\in(\iota,\tau_{2,x})$, pick projections $E_i$ in $\cA_1$ and partial isometries $W_i$ in $\cA_2$ such 
that $\sum_iW_i\circ W_i^*=1_x$ and $W_i^*\circ W_i=\eta(E_i)$, then $\eta_*\sum_i(M^*\circ\mu_2(W_i)\otimes\tau_2(W_i^*)\circ T)=(M^*\otimes T)$ and $\eta_*$ is surjective. 
If $X$ is in the kernel of $\eta_*$, pick $T=\tau_1(E_k)\circ T$ and $M=\mu_1(E_k)\circ M$, where $E_k$ is a minimal projection in $\cA_1$. Then $\eta(E_k)$ is a minimal projection in $\cA_2$ and $T=\tau_2\eta(E_k)\circ T$ and $M=\mu_2\eta(E_k)\circ M$. Thus $\omega_{T,M}(X)=\omega_{T,M}(\eta_*(X)=0$, so $X=0$.\medskip

As a second example of functorial properties we suppose that $\xi:(\sigma,\tilde\sigma)\to(\tau,\tilde\tau)$ 
is a unitary tensor natural transformation. Thus $\xi_\iota=1_\iota$, $\xi_v\circ\sigma(A)=\tau(A)\circ\xi_u$ for $A\in(u,v)$ and 
$\xi_{u\otimes v}\circ\tilde\sigma_{u,v}=\tilde\tau_{u,v}\circ\xi_u\otimes\xi_v$. We now set 
$$\xi_*(M\otimes S):=M\otimes\xi_u\circ S,\quad M\in(\mu_u,\iota),\,\,S\in(\iota,\sigma_u).$$ 
Obviously, if $A\in(u,v)$ then $\xi_*(M\otimes\sigma(A)\circ S)=\xi_*(M\circ\mu(A)\otimes S)$. Thus 
$\xi_*$ can be considered as a linear map from $^\circ_\mu\cC_\sigma$ to  $^\circ_\mu\cC_\tau$. 
$$\xi_*(M\otimes S)\xi_*(M'\otimes S')=(M\otimes M'\circ\tilde\mu_{u,u'}^*\otimes\tilde\tau_{u,u'}\circ(\xi_u\
\circ S)\otimes(\xi_{u'}\circ S')=$$
$$(M\otimes M'\circ\tilde\mu_{u,u'}^*\otimes\xi_{u\otimes u'}\circ\tilde\sigma_{u,u'}\circ S\otimes S')=\xi_*((M\otimes S)(M'\otimes S').$$
Thus $\xi_*$ is multiplicative. 
$$\xi_*((M\otimes S)^*)=\xi_*(M^\bullet\otimes S^\bullet)=(M^\bullet\otimes\xi_{\overline u}\circ S^\bullet)$$ 
whereas 
$$(\xi_*(M\otimes S))^*=(M^\bullet\otimes(\xi_u\circ S)^\bullet).$$ 
Thus it suffices to show that $\xi_u^\bullet=\xi_{\overline u}$. 
$$\xi_{\overline u}^*\circ\xi_u^\bullet=R_{\sigma_u}^*\otimes 1_{\sigma_{\overline u}}\circ 1_{\sigma_{\overline u}}\otimes\xi_u^*\otimes\xi_{\overline u}^*\circ 1_{\sigma_{\overline u}}\otimes\overline R_{\tau_u}.$$
Now 
$$\xi_u^*\otimes\xi_{\overline u}^\circ\tilde\tau_{u,\overline u}*=\tilde\sigma_{u,\overline u}^*\circ\xi_{u\otimes\overline u}^*,\quad \overline R_{\tau_u}=\tilde\tau_{u,\overline u}^*\circ\tau(\overline R_u),$$
and $\xi_{u\otimes\overline u}^*\circ\tau(\overline R_u)=\sigma(\overline R_u)$. This gives $\xi_u^*\otimes\xi_{\overline u}^*=1_{\overline u}$  and hence $\xi_u^\bullet=\xi_{\overline u}$, as required. 
Thus $\xi_*$ is an isomorphism and hence extends to an isomorphism 
from $_\mu\cC_\sigma$ to $_\mu\cC_\tau$.\medskip 

As a third example of functorial properties, we consider quasitensor functors $(\mu,\tilde\mu):\cA\to{\cal M}$, $(\tau,\tilde\tau):\cA\to\cT$ and $(\sigma,\tilde\sigma):\cT\to\cS$ and define for $M\in(\mu_u,\iota)$ and $S\in(\iota,\tau_u)$, $\sigma_*(M\otimes S):=(M\otimes\sigma(S))$. This  obviously defines a linear map $\sigma_*$ from $^\circ_\mu\cC_\tau$ to $^\circ_\mu\cC_{\sigma\tau}$. $\sigma_*$ is multiplicative since 
$$\widetilde{\sigma\circ\tau}_{u,v}\circ\sigma(S)\otimes\sigma(T)=\sigma(\tilde\tau_{u,v})\circ\sigma(S\otimes T),$$ 
where $T\in(\iota,\tau_v)$. Furthermore $\sigma_*$ commutes with the adjoint since $\sigma$ commutes 
with conjugation. Hence $\sigma_*$ extends to a morphism from $_\mu\cC_\tau$ to $_\mu\cC_{\sigma\tau}$.\medskip 

\noindent
{\bf 4.5 Proposition} {\sl $\sigma_*$ is faithful and if $\sigma$ maps 
$(\iota,\tau_u)$ onto $(\iota,\sigma\tau_u)$ for each object $u$ of $\cA$, then $\sigma_*$ is an isomorphism.}\medskip 

\noindent
{\it Proof} Under the above condition, $\sigma_*$ is obviously surjective. It therefore suffices to prove that $\sigma_*$ is faithful on $^\circ_\mu\cC_\tau$. To this end, let us denote the image of an element $X\in\oplus_u(\mu_u,\iota)\otimes(\iota,\tau_u)$ in $\oplus(\mu_u,\iota)\otimes_\cA(\iota,\tau_u)$ 
by $\hat X$ and suppose that $\sigma_*(\hat X)=0$. Then $\sigma_*(X)$, the image of $X$ in $\oplus_u(\mu_u,\iota)\otimes(\iota,\sigma\tau_u)$, is of the form 
$$\sum_i(M_i\otimes (\sigma\tau(A_i)\circ\sigma(S_i))-(M_i\circ\mu(A_i))\otimes\sigma(S_i)).$$
Set $Y:=\oplus_i(M_i\otimes(\tau(A_i)\circ S_i)-(M_i\circ\mu(A_i))\otimes S_i$ then $\sigma_*(X)=\sigma_*(Y)$ and $\hat Y=0$. Hence it suffices to show that $\hat X=0$ when $\sigma_*(X)=0$, but 
this follows since $\sigma$ is faithful on the Hilbert spaces $(\iota,\tau_u)$.\medskip 

We now recall the $^*$--functor $q:\cT\to\cH$, discussed in the appendix, taking an object $x$ of $\cT$ to the Hilbert space $(\iota,x)$ and $X\in(x,y)$ onto the map $T\mapsto X\circ T$. 
$q$ extends uniquely to a quasitensor functor $(q,\tilde q)$,  and $\tilde q$ is minimal. By the above result, $q_*$ is an isomorphism and, in this sense, it 
suffices to consider quasitensor functors $(\mu,\tilde\mu)$ and $(\tau,\tilde\tau)$ taking values in the category of Hilbert spaces.\medskip

\end{section}

\begin{section}{The genesis of ergodic actions} 

In this section we explain how to get actions of quantum groups on the $C^*$--algebras constructed 
from a pair of quasitensor functors in the last section. To this end, we suppose that we have an action $\eta$ of a quantum group $G$ on the category $\cT$ leaving 
 the objects of $\cT$ invariant. Regarding the $C^*$-algebra $\cQ$ of $G$ as a $C^*$--category 
 with a single object, $\eta$ is a $^*$--functor from $\cT$ to $\cT\otimes\cQ$ with $\eta\otimes 1_\cQ\circ\eta=1_\cT\otimes\Delta\circ\eta$, $\Delta$ being the coproduct. Since $\cT\otimes\cQ$ is not a tensor $C^*$--category, we cannot require $\eta$ to be a tensor $^*$--functor. The natural condition is to require 
instead that $\eta(S\otimes T)=\eta(S)\top\eta(T)$, where $\top$ indicates that we take a tensor product in the first component and a product in the second a la Woronowicz. This is the product used when 
defining the tensor product of representations of a quantum group. We further suppose that the arrows 
of $\cA$ intertwine this action, i.e.\ that $\eta(\tau(A)\circ T)=\tau(A)\otimes I\circ\eta(T)$ and 
 that $\eta(\tilde\tau_{u,v})=\tilde\tau_{u,v}\otimes I$ for each pair $u,v$ of objects of $\cA$.\medskip 
 
\noindent
{\bf 5.1 Proposition} {\sl There is a unique action $\alpha$ of $G$ on $_\mu\cC_\tau$ such 
 that $\alpha(M\otimes T):=M\otimes\eta(T)$.} \medskip

 \noindent
 {\it Proof} $\alpha$ is obviously well defined and a simple computation shows that it is multiplicative.
To show that it commutes with the adjoint, we must show that $\eta(S^\bullet)=\eta(S)^{\bullet\otimes*}$. 
Now $$\eta(S^\bullet)=\eta(S^*\otimes 1_{\tau_{\overline u}}\circ\tilde{\overline R}_u)=\eta(S^*)\top
1_{\tau_{\overline u}}\otimes I\circ\tilde{\overline R}_u\otimes I=\eta(S)^{\bullet\otimes*}.$$ 
Finally, $\alpha$ extends to $_\mu\cC_\tau$ by continuity and is trivially an action.\medskip  

Let $E_G$ denote the conditional expectation defined by averaging $\eta$ over $G$ and let $E_{G,u}$ 
denote the projection obtained by restricting to $(\iota,\tau_u)$. If  $A\in(u,v)$ then $E_{G,v}\tau(A)E_{G,u}=\tau(A)E_{G,u}$ thus $E_{G,v}\tau(A)=\tau(A)E_{G,u}$. Thus $E_{G,\cdot}$ is a natural transformation from $\tau$ to $\tau$ as is $\tau(c)$ where $c$ denote the central support of $\iota$ in $\cA$ as before. If $E_{G,\cdot}=\tau(c)$ then the above action will be ergodic.\medskip

It remains to understand how to get appropriate actions of a quantum group $G$ on $\cT$ 
and when the induced action on $_\mu\cC_\tau$ is ergodic. Now if we suppose that $\tau$ 
is a tensor $^*$--functor into the category of Hilbert space then, as $\cA$ has conjugates, 
the duality theorem of Woronowicz \cite{WoronowiczTK} gives us an action $\eta$ of a compact quantum group $G_\tau$ on the Hilbert spaces $(\iota,\tau_u)$ and hence on the category $\cT$. The $C^*$--algebra of $G_\tau$ is $_\tau\cC_\tau$, showing how our construction generalizes that of Woronowicz. 
The action is defined by 
$$\eta_u(T)=\sum_iT_i\otimes(T_i^*\otimes T),\quad T\in(\iota,\tau_u),$$ 
where the sum is taken over an orthonormal basis. The arrows of the form $\tau(A)$ intertwine this action, the conditional expectation defined by averaging over $G_\mu$ is $\tau(c)$ and $\eta(S\otimes T)=\eta(S)\top\eta(T)$. Thus we have an induced action $\alpha$ of $G_\tau$ on $_\mu\cC_\tau$: 
$$\alpha(M\otimes T)=\sum_i(M\otimes T_i)\otimes(T_i^*\otimes T).$$ 
This action is ergodic since the conditional expection coincides with $h$ which is 
the unique  invariant state. To put ourselves in the setting of \cite{PR}, we 
replace $\mu$ by $q\mu$ which can then be identified with the spectral 
functor of the ergodic action. The spectral space $\overline{L_u}$ associated with
the representation $u$ can be identified, as a Hilbert space, with
$(\iota,\mu_u)$ through the map that takes $M\in(\iota,\mu_u)$
to the complex conjugate of the map of $L_u$:
$$T\in(\iota,\tau_u)\to M^*\otimes T\in{}_\mu\cC_\tau.$$

\noindent
{\it Remark} Under certain circumstances, the above generalizes to the case of a quasitensor functor $\tau$ but since it is not needed here, we will give details in a separate paper.\medskip

\end{section}

\begin{section}{Self-conjugate solutions of the conjugate equations} 

To treat self-conjugate solutions of the conjugate equations we
consider two tensor $^*$--categories. The first, $\cT_{rd}$, for real solutions, has objects that are powers $y^n$, $n\in\mathbb N_0$, of a generating object $y$ and whose arrows are generated by a single arrow $S\in(\iota,y^2)$ satisfying $S^*\otimes 1_y\circ 1_y\otimes S=1_y$ and $S^*\circ S=d$. The second, $\cT_{pd}$, for pseudoreal solutions, has objects that are powers $z^n$, $n\in\mathbb N_0$, of a generating object $z$ and whose arrows are generated by a single arrow $S$ satisfying $S^*\otimes 1_y\circ  1_y\otimes S=-1_y$ and $S^*\circ S=d$. In both cases we suppose that $d\neq 0$.\medskip 

As $\cT_{rd}$ and $\cT_{pd}$ are defined in terms of generators and relations  they will satisfy the 
corresponding universal properties. But the analogous universal properties are satisfied by the 
Temperley-Lieb categories \cite{Y}, usually defined without reference to a $^*$--operation. Hence 
the categories $\cT_{rd}$ and $\cT_{pd}$ are Temperley-Lieb categories corresponding to parameters $\pm d$ with a $^*$--operation defined by a solution $S$ of  self-conjugate solutions 
of the conjugate equations. As such, the following assertions  are well known. The units and generating objects of these categories are irreducible and the spaces of arrows are finite dimensional. The categories are simple except at roots of unity, $d=2\cos{\pi\over\ell}$, $\ell=3,4,\dots$, when they have a single non-zero proper ideal \cite{GW}. They are tensor $C^*$--categories 
when $d\geq 2$ and at roots of unity their quotients by the unique non-zero proper ideal are tensor 
$C^*$--categories having the universal property, but now for tensor $C^*$--categories.\medskip 

We define a left inverse $\psi$ for the generating object $y$ of $\cT_{rd}$ by 
$$\psi_{m,n}(Y):=S^*\otimes 1_{y^{m-1}}\circ 1_y\otimes Y\circ 1_{y^{n-1}}\otimes S,\quad Y\in(y^n,y^m).$$ 
Iterating $\psi$ we get a mapping Tr$:(y^n,y^n)\to(\iota,\iota)$, the Markov 
trace. A right inverse for $y$ is obtained by dualizing the above definition 
with respect to $\otimes$ and iterating again defines the Markov trace.\medskip 

A left and right inverse for the generating object $z$ of $\cT_{pd}$ can be 
defined analogously and their iterates again yield the Markov trace.\medskip

\end{section} 

\begin{section}{General solutions of the conjugate equations}

 We define $\cT_d$  for $d\neq 0$ to be the tensor $^*$--category whose objects are the words                                                                                                                                                                                                                                                                                                                                                                                                                                                                                                                                                                                                                                                                          
in $x$ and $\overline x$  and whose 
arrows are generated by two arrows $R\in(\iota,\overline x\otimes x)$ and 
$\overline R\in(\iota,x\otimes \overline x)$ subject to the relations 
$1_x\otimes R^*\circ \overline R\otimes 1_x=1_x$, 
$1_{\overline x}\otimes\overline R^*\circ R\otimes 
1_{\overline x}=1_{\overline x}$, $R^*\circ R=d$ and 
$\overline R^*\circ\overline R=d$.
We note that $\cT_d$ has an involution $^\circ$ taking $x$ to $\overline x$ and $R$ to 
$\overline R$.
We call an object even or odd according as to whether it a tensor product of an even or odd number of the objects $x$ and $\overline x$. The space of arrows between an even and an odd object is zero.\medskip 

As $\cT_d$ has been defined in terms of generators and relations, it has a universal property: 
given any solutions $R',\overline R'$ of the conjugate equations in a tensor $^*$--category $\cT$ with 
$R^{'*}\circ R'=d$ and $\overline R^{'*}\circ\overline R'=d$, there is a unique tensor $^*$--functor  $\phi:\cT_d\to\cT$ such that $\phi(R)=R'$ and $\phi(\overline R)=\overline R'$. Yamagami in \cite{Y1} defines a tensor 
$^*$--category in terms of oriented Kauffman diagrams and shows that it has the above universal 
property so that this category is in fact isomorphic to $\cT_d$. With his very different starting point his 
proof of  Theorem 7.1 below is quite different.\medskip 

\noindent
{\bf 7.1 Theorem} {\sl $\iota$ and $x$ are irreducible in $\cT_d$.}\medskip 

\noindent
{\it Proof} If $X\in(\iota,\iota)$ then $X\to 1_x\otimes X$ is a faithful 
morphism from $(\iota,\iota)$ to $(x,x)$, since it suffices to show that $(\iota,x\otimes\overline x)$ is $1$-dimensional in $\cT_d$. This will be 
  the case if any intertwiner in $(\iota,x\otimes\overline x)$ constructed as an algebraic 
  expression in $R$, $\overline R$, $1_x$, $1_{\overline x}$ and their adjoints reduces 
  to a multiple of $\overline R$. A tensor product of the basic arrows $1_x$, 
  $1_{\overline x}$, $R$, $\overline R$, $R^*$ and $\overline R^*$. 
  will be said to be a term. A term is positive if $R^*$ and $\overline R^*$ are not 
  involved and negative if $R$ and $\overline R$ are not involved. Using the interchange law, 
  any term can be written in the form $X_+\circ X_-$, where $X_+$ is a positive term 
  and $X_-$ a negative term. Now consider a composition of terms of the form $X_-\circ X_+$. 
   We break these two terms into an equal number of pieces of minimal size such that the 
  $\circ$--composition of the corresponding pieces is defined. We list the possible 
  $\circ$--compositions of two pieces. 
  $1_ x\otimes R^*\circ \overline R\otimes 1_x=1_x$; 
  $R^*\otimes 1_{\overline x}\circ 1_{\overline x}\otimes\overline R=1_{\overline x}
  $; $1_{\overline x}\otimes\overline R^*\circ R\otimes 1_{\overline x}
  =1_{\overline x}$; $\overline R^*\otimes 1_x\circ 1_x\otimes R=1_x$; 
  $R^*\circ R=d$; $\overline R^*\circ\overline R=d $; 
    $1_x\circ 1_x=1_x$; $1_{\overline x}\circ 1_{\overline x}=1_   {\overline x}$; $1_{\overline x\otimes x}\circ R=R$; $1_{x\otimes\overline x}\circ \overline R=\overline R$; $R^*\circ 1_{\overline x\otimes x}=R^*$; $\overline R^*\circ 1_{x\otimes\overline x}=\overline R^*$; 
 $1_{x\otimes\overline x}\otimes R^*\circ\overline R\otimes 1_{\overline x\otimes x}=R\otimes R^*$; 
  $1_{\overline x\otimes x}\otimes R^*\circ R\otimes 1_{\overline x\otimes x}=R\otimes R^*$; 
  $1_{x\otimes\overline x}\otimes \overline R^*\circ\overline R\otimes 1_{x\otimes\overline x}=\overline R\otimes \overline R^*$; $1_{\overline x\otimes x}\otimes\overline R^*\circ R\otimes 1_{x\otimes\overline x}=R\otimes R^*$. 
Thus up to a scalar $X_-\circ X_+=Y_+\circ Y_-$. Hence up to a scalar any composition of terms 
can be written as a composition of positive and negative terms where the negative terms appear 
on the right and the positive terms on the left. But such a composition is an arrow of $(\iota,x\otimes\overline x)$ 
if and only if there are no negative terms and a single positive term $\overline R$.\medskip

\noindent
{\it Remark} Every composition of terms in $\cT_d$ can also be written as a composition of positive and negative terms, the positive terms appearing on the right and the negative terms on the left.\medskip

The universal property of $\cT_d$ implies that there is a unique tensor $^*$-functor $\phi:\cT_d\to\cT_{rd}$ such that $\phi_x=y$ and $\phi(R)=\phi(\overline R)=S$. As an aid to studying this functor we 
introduce the full subcategory $\cT_d^a$ whose objects are of the form $(x\otimes\overline x)^n$, the even objects or $(x\otimes\overline x)^n\otimes x$, the odd objects, for $n\in\mathbb N_0$. This category is obviously not a tensor 
subcategory of $\cT_d$ but it can be given the structure of a tensor $^*$--category. Note first that 
there is no non-zero arrow between objects of different parity so any non-zero arrow has a definite 
parity and we define the tensor product as $X\otimes X'$ if $X$ is even and $X\otimes X^{'\circ}$ if 
$X$ is odd. As an object of $\cT_d^a$, $x$ is self-conjugate. Therefore there is a  unique tensor $^*$--functor $\psi:\cT_{rd}$ to $\cT_d^a$ with $\psi_y=x$ and $\psi(S)=\overline R$. The restriction $\phi^a$ of $\phi$ to $\cT_d^a$ is a tensor $^*$--functor since $\phi(X)=\phi(X^\circ)$ and $\phi^a\psi=1_{\cT_{rd}}$. Thus $\phi^a$ and $\phi$ are surjective on arrows. $\phi$ is obviously not full since $(x,\overline x)=0$ and $(y,y)\neq 0$. Thus, in particular,  we have proved the following result.
\medskip

\noindent
{\bf 7.2 Theorem} {\sl The canonical functor $\phi$ from $\cT_d$ to $\cT_{rd}$ is a tensor $^*$--functor surjective on both objects and arrows but not full.}\medskip

\noindent
{\bf 7.3 Theorem} {\sl $\psi:\cT_{rd}\to \cT_d^a$ is an isomorphism of tensor $^*$--categories and 
$\phi:\cT_d\to\cT_{rd}$ is faithful.}\medskip 

\noindent
{\it Proof} If $X$ is a term of $\cT_d^a$ then $\phi^a(X)$ will be a term of $\cT_{rd}$ and $\psi\phi^a(X)=X$. But every arrow of $\cT_d^a$ is a linear combination of compositions of terms. Hence $\psi$ is 
an isomorphism. More generally, given a full subcategory $\cT_d^s$ of $\cT_d$ such that the restriction $\phi^s$ of $\phi$ to the objects of $\cT_d^s$ is an isomorphism, the image of $\cT_d^s$ under $\phi^s$ is a subcategory $\cT_{rd}^s$ of $\cT_{rd}$. Every term $Y$ of $\cT_{rd}^s$ is the image under $\phi^s$ of a term of unique term $\psi^s(Y)$ of $\cT_d^s$. Note, however, that formally distinct terms of $\cT_d$ can define the same arrow of $\cT_d$, for example, $R\circ\overline R^*=\overline R^*\otimes R$. Nevertheless $\psi^s$ extends to a full functor from $\cT_{rd}^s$ to $\cT_d^s$. $\phi^s$ and $\psi^s$ are isomorphisms since they are inverses of one another in restriction to terms. Using the linear isomorphism 
of $(x\otimes p,q)$ and $(p,\overline x\otimes q)$ where necessary, we conclude that $\phi$ is faithful.\medskip

We define left inverses $\psi$ and $\overline\psi$ of $x$ and $\overline x$ by 
$$\psi_{p,q}(X):=R^*\otimes 1_p\circ 1_{\overline x}\otimes X\circ R\otimes 1_q,\quad X\in(x\otimes q,x\otimes p),$$ 
$$\overline\psi_{p,q}(X):\overline R^*\otimes 1_p\circ 1_x\otimes X\circ\overline R\otimes 1_q.\quad X\in(\overline x\otimes q,\overline x\otimes p).$$ 
Iterating $\psi$ and $\overline\psi$ appropriately we get a map Tr$:(p,p)\to(\iota,\iota)$. Obviously, 
since the unit of $\cT_d$ is irreducible by Theorem 7.1, Tr$(X)=$Tr$(\phi(X))$. Right inverses of $x$ and $\overline x$ are defined by dualizing with respect to $\otimes$ and iterated appropriately lead to the same map Tr$:(p,p)\to(\iota,\iota)$. It follows that Tr is a trace in that Tr$(X\circ X')=$Tr$(X'\circ X)$, 
whenever the compositions are defined. Another way of looking at this trace is to note that $R,\overline R$ extend uniquely to a homomorphic choice $q\mapsto R_q$ of solutions of the conjugate equations in $\cT_d$ (see Appendix). The corresponding scalar product $\psi_q$ on $(q,q)$ is that associated with the trace.\medskip 

Following \cite{Wenzl3}, we call an arrow $X\in(p,q)$ of $\cT_d$ negligible if Tr$(X'\circ X)=0$ for 
all $X'\in(q,p)$. Clearly, $X$ is negligible if $\phi(X)$ is negligible in $\cT_{rd}$ but the converse follows 
from Theorem 7.2. The set of negligible arrows is a tensor $^*$--ideal. This is well known in the 
case of $\cT_{rd}$ but holds in some generality.\medskip 

\noindent
{\bf 7.4 Proposition} {\sl Let $\cT$ be a tensor $^*$--category with conjugates and irreducible tensor unit and $u\mapsto R_u$ a tracial and homomorphic choice of solutions of the conjugate equations with   associated trace Tr, then the set $\cI$ of negligible arrows in $\cT$ is the maximal proper tensor $^*$--ideal.}\medskip 

\noindent
{\it Proof} If $X\in(u,v)$ then Tr$(1_u)=R_u^*\circ 1_{\overline u\otimes u}\circ R_u=d_u
$. So Tr$(X'\circ X)^*=$Tr$X^*\circ X^{'*}=$Tr$X^{'*}\circ X^*$ so $\cI=\cI^*$. If $W\in(w,u)$ and $V\in(v,w)$ then Tr$(V\circ X\circ W)=$Tr$(W\circ V\circ X)$ so $X\in\cI$ implies $X\circ W\in\cI$. Now let $Z\in(t\otimes v,t\otimes u)$, then, since $u\mapsto R_u$ 
is homomorphic,
$${\rm Tr}(Z\circ 1_t\otimes X)=R_u^*\circ 1_{\overline u}\otimes(R^*_t\otimes 1_t\circ 1_{\overline t}\otimes Z\circ R_t\otimes 1_v\circ X)\circ R_u.$$ 
Thus $X\in\cI$ implies $1_t\otimes X\in\cI$ and similarly for tensoring on the right. Tr$(1_\iota)=1$ so that $\cI$ is a proper tensor $^*$--ideal. Now if $\cJ$ is a proper tensor $^*$-ideal and $X\in(u,v)\cap\cJ$ then Tr$(X'\circ X)\in\cJ$ for all $X'\in(v,u)$ thus $X\in\cI$, completing the proof.\medskip 

\noindent
{\bf 7.5 Proposition} {\sl A tensor $C^*$--category with irreducible tensor unit and conjugates is simple.}\medskip

\noindent
{\it Proof} A standard choice of solutions of the conjugate equations yields a trace independent of that 
choice (see Appendix) so the proof of  Proposition 7.4 applies. Since the trace is faithful every 
negligible arrow is zero.\medskip 

We now investigate the ideal structure of $\cT_d$. As remarked above, an arrow $X$ of $\cT_d$ is 
negligible if and only if $\phi(X)$ is negligible. Thus if $d\neq\cos{\pi\over\ell}$ for $\ell=3,4,\dots$, 
$X$ is negligible if and only if $\phi(X)=0$. Thus by Theorem 7.3, $\cT_d$ is simple for these  values of $d$. If $d=\cos{\pi\over\ell}$ for $\ell=3,4,\dots$, then again by Theorem 7.2, the ideal of negligible arrows in $\cT_d$ is the unique non-zero proper tensor ideal. We have therefore proved the following result.\medskip 

\noindent
{\bf 7.6 Theorem} {\sl $\cT_d$ is simple unless $d=2\cos{\pi\over\ell}$ when it has a single non-zero 
proper tensor ideal.}\medskip

As another consequence of Theorem 7.3, $\cT_d$ will be a tensor $C^*$--category whenever $\cT_{rd}$ 
is a tensor $C^*$--category, i.e.\ whenever $d\geq 2$. If $d=2\cos{\pi\over\ell}$, $\ell=3,4,\dots$, its quotient by the unique non-zero proper ideal will be a tensor $C^*$--category having the universal 
property for normalized solutions of the conjugate equations with these values of $d$.\medskip 

There is also a canonical functor $\phi:\cT_d\to\cT_{pd}$ with $\phi_x=\phi_{\overline x}=z$, 
$\phi(R)=S$ and $\phi(\overline R)=-S$. The results of this section have obvious analogues in 
this case.\medskip

\end{section}

\section{Embedding the universal categories} 

There will be tensor $^*$--functors from $\cT_d$ to a tensor $C^*$--category $\cT$ if $d=1$, $d\geq 2$ or if $d=2\cos{\pi\over\ell}$, for $\ell=3,4,\dots$. For these values of $d$, 
given normalized solutions $R',\overline R'$ of the conjugate equations for $x'$ in $\cT$ with $R^{'*}\circ R'=d$, there is a unique tensor $^*$--functor $\phi$ from $\cT_d$ to $\cT$ with $\phi_x=x'$, $\phi(R)=R'$ and $\phi(\overline R)=\overline R'$. When $\phi$ is injective on objects the image of $\phi$ is a tensor $C^*$--category isomorphic to $\cT_d$. 
Similarly, if $S'\in(\iota,y^{'2})$ is a real (or pseudoreal) solution of the conjugate
equations in ${\cal T}$ with $S'^*\circ S'=d$, then there is a unique isomorphism $\phi$ from $\cT_{rd}$ (or $\cT_{pd}$) to the tensor $C^*$--subcategory
generated by $S$ taking $y$ to $y'$ and $S$ to $S'$.\medskip

Obviously, a tensor $C^*$--category $\cT$ with $(\iota,\iota)=\mathbb C$ can only be 
embedable in the tensor $C^*$--category of Hilbert spaces if  any normalized 
solution $R,\overline R$ of the conjugate equations in $\cT$ with $R^*\circ R<2$ is 
unitary. In particular, the quotient tensor $C^*$--categories $\cT_d/\cI$ for $d=2\cos{\pi\over\ell}$ cannot 
be embedded into the tensor $C^*$--category of Hilbert spaces.\medskip

We will now classify, for the possible values of $d$ other than $d=1$, the tensor $^*$--functors from  $\cT_{rd}$ and $\cT_{pd}$ to the category of Hilbert spaces up to natural unitary tensor equivalence. They are 
determined by the parameters $\lambda_i$ used by \cite{BRV}. Let $\phi:\cT_{rd}\to\cH$ be a tensor $^*$--functor, then there are $0<\lambda_i<1$ with $\sum_1^k(\lambda_i^2+\lambda^{-2}_i)+n-2k=d$ and an orthonormal basis $e_i$ of $\phi_\sigma$ such that 
$$\phi(S)=\sum_1^k\lambda_ie_{i+k}\otimes e_i+\sum_1^k\lambda^{-1}_ie_i\otimes e_{i+k}+\sum_{2k+1}^ne_i\otimes e_i.$$ 
Similarly, if $\phi:\cT_{pd}\to\cH$ is a tensor $^*$--functor then there are $0<\lambda_i\leq 1$ with $\sum_1^{n\over 2}(\lambda_i^2+\lambda^{-2}_i)=d$ and an orthonormal basis $e_i$ of $\phi_\sigma$ such that 
$$\phi(S)=\sum_1^{n\over 2}\lambda_ie_{i+{n\over 2}}\otimes e_i-\sum_1^{n\over 2}\lambda_i^{-1}e_i\otimes e_{i+{n\over 2}}.$$
The following result can easily be proved, cf. \cite{Y2}.\medskip

\noindent 
{\bf 8.1 Proposition} {\sl 
\begin{description}
\item{a)}The parameters $\lambda_i$ with $0<\lambda_i<1$ and  $\sum_1^k(\lambda_i^2+\lambda^{-2}_i)+n-2k=d$ classify the tensor $^*$--functors $\phi:\cT_{rd}\to\cH$ 
up to a natural unitary tensor equivalence.
\item{b)}The parameters $\lambda_i$ with $0<\lambda_i\leq 1$ and  $\sum_1^{n\over 2}(\lambda_i^2+\lambda^{-2}_i)=d$ classify the tensor $^*$--functors $\phi:\cT_{pd}\to\cH$ up to a tensor unitary natural equivalence. 
\end{description} }\medskip 

In case a), the spectrum of $j^*_y\circ j_y$ is $\{\lambda^2_i,\lambda^{-2}_i:1\leq i\leq k\}\cup\{1:2k+1\leq i\leq n\}$, where $j_y$ is the antilinear invertible operator  on $H_y$ defined by $\phi(S):=\sum_i e_i\otimes  j_ye_i$.  In case b) it is $\{\lambda^2_i,\lambda^{-2}_i:1\leq i\leq{n\over 2}\}$.\medskip

We next classify the tensor $^*$--functors $\phi:\cT_d\to\cH$, $d\neq 1$, up to a tensor unitary natural equivalence. Any such functor $\phi$ determines an invertible antilinear operator $j_x$ on $H_x$ via 
$\phi(\overline R_x):=\sum_i e_i\otimes j_x e_i$. The following result can again be easily proved.
\medskip
 
\noindent
{\bf 8.2 Proposition} {\sl The tensor natural unitary equivalence 
classes of embeddings of $\cT_d$ 
into Hilbert spaces are classified by a monotone set of parameters $0<\lambda_i$ with $\sum_{i=1}^n\lambda_i^2=\sum_{i=1}^n\lambda_i^{-2}=d$, where $\{\lambda_i^2\}$ is just the eigenvalue 
list of $j_x^*\circ j_x$.}\medskip 

A rather less natural description of these equivalence classes in terms of an invertible linear operator 
can be found in \cite{Y1}.\medskip 

\noindent
{\it Remark} As we have canonical tensor $^*$--functors from $\cT_d$ to $\cT_{rd}$ and $\cT_{pd}$ 
taking $R$  and $\overline R$ onto $S$ and $-S$ respectively, an embedding of $\cT_{rd}$ or $\cT_{pd}$ induces an 
embedding of $\cT_d$, equivalent embeddings inducing equivalent embeddings. Thus for each set of parameters in Proposition 8.1 there is a corresponding set of parameters in Proposition 8.2. The eigenvalue list of $j_y^*j_y$ in Proposition 8.1a and $j_z^*j_z$  in Proposition 8.1b has been indicated above. 
Recall that $n$ is even in Proposition 8.1b. We see that no two inequivalent embeddings of $\cT_{rd}$ 
or $\cT_{pd}$ can induce equivalent embeddings of $\cT_d$ but that for each embedding of 
$\cT_{pd}$ there is an embedding of $\cT_{rd}$ inducing an equivalent embedding of $\cT_d$.\medskip

  If $\cA$ is a tensor $C^*$--category with conjugates and irreducible tensor unit, then by the duality theorem of Woronowicz \cite{WoronowiczTK} every embedding of $\cA$ into the category of Hilbert spaces determines a compact 
  quantum group whose category of finite-dimensional representations is equivalent to the completion 
  of $\cA$ under subobjects and direct sums. It can be easily shown that embeddings differing by a 
  tensor unitary equivalence yield isomorphic quantum groups. In fact, 
  this results from the discussion following Proposition 4.4.\medskip 
 
   We now describe the compact quantum groups arising from the embeddings of $\cT_d$, $\cT_{rd}$ and $\cT_{pd}$. Obviously all of these quantum groups have a representation theory  generated by a single fundamental representation $u$. In the Hilbert space of the embedding, the arrows $R\in(\iota,\overline u\otimes u)$ and $\overline R\in(\iota,u\otimes\overline u)$ are then intertwining operators between the 
  associated representations of the compact quantum group. The conjugation $j_u:H_u\to H_{\overline u}$
is an invertible antilinear intertwiner $j_u$ from $u$ to $\overline u$. Thus $j_u\otimes^*\circ u=\overline u\circ  j_u$, or, in terms of matrix elements $j_{mn}u_{np}^*=\overline u_{mr}j_{rp}$. This relation might be used to define the compact quantum groups involved. In fact,  the quantum groups involved have been defined, less intrinsically, in terms of a linear operator $Q$ in the notation of Wang \cite{Wang2} and $F$ in that of Banica \cite{Banica0}, \cite{Banica}.\medskip 
 
 We first treat the self-conjugate case, i.e.\ embeddings of $\cT_{rd}$ and $\cT_{pd}$ so that
 $\overline u=u$ and hence $j_u^2=\pm 1$. We let $c$ be an antiunitary involution  on the Hilbert  space of $u$ and set $Q:=cj_u^*$, then $QcQc=\pm 1$ and $Q^*\otimes 1\circ c\otimes^*\circ u=u\circ Q^*\circ c$, then working in the basis where $c_{ij}=\delta_{ik}$, we get 
$u_{np}Q_{nm}=Q_{pn}u_{mn}^*$ or $u^t\circ Q=Q\otimes 1\circ u^*$.  These are the defining relations for the compact quantum group $B_u(Q)$ in the notation of Wang. Note that $Q^*Q=j_uj_u^*$, thus the isomorphism class of $B_u(Q)$ depends only on the eigenvalue list of $Q^*Q$, improving Wang's result. Banica uses the adjoint operator $F:=Q^*=j_uc$ and denotes the quantum group by $A_o(F)$.

Turning to the embeddings of $\cT_d$, Banica and Wang make use not of the conjugate representation $\overline u$ but of the equivalent non-unitary representation $\tilde u$, $\tilde u_{mn}:=u_{mn}^*$ which depends on a choice of orthonormal basis in $H_u$. The antiunitary operator $c$ leaving this 
basis fixed intertwines $\tilde u$ and $u$. Thus $F:=j_uc$ is a linear intertwiner from $\tilde u$ and $\overline u$ and, setting $Q:=F^*F$, a computation shows that $u_{np}Q_{nr}=Q_{ps}u_{rs}^*$ or 
$u^t\circ Q=Q\otimes 1\circ\tilde u$. This is the relation used by Wang to define the compact quantum group $A_u(Q)$ or $A_u(F)$ in the notation of Banica. Note that the eigenvalue list of $Q$ coincides with that of $j_u^*j_u$ and is hence 
characteristic of the tensor natural unitary equivalence class of the embedding. As Wang showed, the quantum groups $A_u(Q)$ and $A_u(Q^{-1})$ are isomorphic and this reflect the involution on $\cT_d$ exchanging $R$ and $\overline R$. The relation between the groups $A_u(Q)$  and the embeddings 
of $\cT_d$ have already been established by Yamagami \cite{Y1}.\medskip

  Thus given a normalized solution of the conjugate equations $R,\overline R$ in a tensor $C^*$--category $\cM$ we have a canonical tensor $^*$--functor $\mu:\cT_d\to\cM$ and picking 
  an embedding $\tau$ into the category of Hilbert spaces we get an ergodic action of $G_\tau$ on $_\mu\cC_\tau$. Choosing $\tau$ suitably, $G_\tau\simeq A_u(Q)$ for any $Q>0$ with Tr$(Q)=R^*\circ R$.\medskip 
  
  Similarly, given a real solution of the conjugate equations $R$ in a tensor $C^*$--category $\cM$, we have a canonical tensor $^*$--functor $\mu:\cT_{rd}\to\cM$ and picking an embedding $\tau$ of $\cT_{rd}$ into the category of Hilbert spaces, we get an ergodic action of $G_\tau$ on $_\mu\cC_\tau$. 
Choosing $\tau$ suitably, $G_\tau\simeq B_u(Q)$ for any $Q$ with Tr$(Q^*Q)={\rm Tr}(Q^*Q)^{-1}=R^*\circ R$ and 
$QcQc=I$. Since $R=\sum_k\psi_k\otimes j_u\psi_k$, where the sum is taken over an orthonormal basis of $\tau_u$  invariant under $c$, we have $R=\sum_k\psi_k\otimes Q^*\psi_k$.\medskip 

  Given a pseudoreal solution of the conjugate equations $R$ in a tensor $C^*$--category $\cM$, we have a canonical tensor $^*$--functor $\tau:\cT_{pd}\to\cM$ and picking an embedding $\tau$ of 
  $\cT_{pd}$ into the category of Hilbert spaces, we get an ergodic action of $G_\tau$ on $_\mu\cC_\tau$. Choosing $\tau$ suitably, $G_\tau\simeq B_u(Q)$ for any $Q$ with Tr$(Q^*Q)={\rm Tr}(Q^*Q)^{-1}=R^*\circ R$ and 
  $QcQc=-I$. The comment on ergodic actions of S$_\mu$U(2) follows since S$_\mu$U(2) is isomorphic to $B_u(Q)$ with $QcQc=1$ when $\mu>0$ and the eigenvalue list of $Q^*Q$ is $|\mu|\leq|\mu^{-1}|$ and to $B_u(Q)$ with $QcQc=-1$ when $\mu<0$ and the eigenvalue list of $Q^*Q$ is $\mu<\mu^{-1}$. In this case, $R=-\sum_k\psi_k\otimes Q^*\psi_k$.\medskip 
  
 If we pick $v\mapsto R_v$ to be standard,then the condition $m(v)={\rm dim}(\iota,\mu_v)$ is equivalent  to saying that $\hat R_v$ is standard.   The results on the q-multiplicity now follow from Corollary A.10.
This completes the proof of Theorems 3.1 and 3.2.\medskip

\begin{section}{Outlook} 

The work reported on in this paper is in the process of being extended in several directions.
In Sections 6, we introduced the tensor $^*$--categories with conjugates $\cT_{rd}$ and $\cT_{pd}$ 
whose objects were tensor powers of a single object and described their embeddings into 
the tensor category of Hilbert spaces and the associated compact quantum groups. An interesting 
problem is to describe tensor $C^*$--categories without conjugates whose objects are again tensor powers of a single irreducible object but where the completion under subobjects has conjugates, their 
embeddings into Hilbert spaces and the associated compact quantum groups. The compact quantum 
groups S$_\mu$U(n), $n\geq 3$ are the prime examples that can be obtained in this way. We have not found a systematic way of producing further examples nor of classifying the underlying tensor $C^*$--categories. However, we have found compact quantum groups depending on two integers $n>2$, the 
smallest integer $n>0$ such  that $\iota\leq x^n$, where $x$ is the generating object  and $d$ is the 
intrinsic dimension of $x$ and also the dimension of the Hilbert space of the corresponding 
representation of the compact quantum group.\medskip 

An interesting aspect, not touched on in the paper, of our way of constructing  ergodic actions, is 
that it really leads to two ergodic actions on the $C^*$--algebra ${}_\mu\cC_\tau$. Thus if $G_\mu$ and $G_\tau$ denote the quantum groups with $\cC(G_\mu)={}_\mu\cC_\mu$ and $\cC(G_\tau)={}_\tau\cC_\tau$ then $G_\mu$ acts on the left and $G_\tau$ on the right on $_\mu\cC_\tau$. The simplest and well known example of this phenomenon is when $\mu=\tau$ yielding the left and right actions of a 
quantum group on itself.\medskip 

Finally, the $C^*$-algebras ${}_\mu\cC_\tau$ may be used in a different way; we may define a 
suitable left action of the algebra on itself together with the obvious right action, making it into
a ${}_\mu\cC_\tau$--bimodule. Further bimodules can be constructed too, and 
reflect more fully the structure of the target tensor $C^*$--category $\cM$ of  $\mu$.  When $\mu$ is 
surjective on objects, $\cM$ can be embedded in a tensor $C^*$--category of ${}_\mu\cC_\tau$--bimodules and this in turn leads 
to further insight on when $\cM$ can be embedded into the tensor $C^*$--category of Hilbert spaces.\medskip

\end{section}

\begin{section} {Appendix}

In this section we show some properties of quasitensor functors and conjugation used in the paper 
and begin by establishing the equivalence of the definition of quasitensor functor with that in \cite{PR}.

Composing $(2.9)$ on the left by $1_{\tau_u}\otimes\tilde\tau_{v,w}^*$ and on the right by $\tilde\tau_{u,v}\otimes 1_{\tau_w}$ gives
$$1_{\tau_u}\otimes\tilde\tau_{v,w}^*\circ\tilde\tau_{u,v\otimes w}^*\circ\tilde\tau_{u\otimes v,w}\circ \tilde\tau_{u,v}\otimes 1_{\tau_w}=1_{\tau_u}\otimes 1_{\tau_v}\otimes 1_{\tau_w}.$$
Since we are dealing with isometries, this implies
$$\tilde\tau_{u\otimes v, w}\circ \tilde\tau_{u,v}\otimes 
1_{\tau_w}=\tilde\tau_{u,v\otimes w}\circ 1_{\tau_u}\otimes 
\tilde\tau_{v,w}=:\tilde\tau_{u,v,w},\eqno(A.1)$$ 
the associativity condition.
If we let $E_{u,v}\in(\tau_{u\otimes v}, \tau_{u\otimes v})$ be the range 
projection  of
$\tilde\tau_{u,v}$ and 
$E_{u,v,w}\in(\tau_{u\otimes v\otimes 
w},\tau_{u\otimes v\otimes w})$
be the range projection of $\tilde\tau_{u,v,w}$, 
then by (2.10) and (A.1)
$$E_{u,v\otimes w}\circ E_{u\otimes v,w}=\tilde\tau_{u,v\otimes w}\circ 1_{\tau_u}\otimes\tilde\tau_{v,w}
\circ\tilde\tau_{u,v}^*\otimes 1_{\tau_w}\circ\tilde\tau_{u\otimes v,w}^*$$
$$=\tilde\tau_{u\otimes v,w}\circ\tilde\tau_{u,v}\otimes 1_{\tau_w}\circ{\tilde\tau_{u,v}}^*\otimes 1_{\tau_w}\circ\tilde\tau_{u\otimes v,w}^*=E_{u,v,w}.\eqno(A.2)$$ 
(A.1) and (A.2) replaced (2.10) in the definition of quasitensor functor in \cite{PR}. 
On the other hand, composing (A.2) on the left with $\tilde\tau_{u\otimes v,w}^*$ and on 
the right with $\tau_{u,v\otimes w}$ and using (A.1), we get (2.10). Thus the two definitions 
are equivalent.\medskip

\noindent
{\it Remark} We automatically have $\tilde\tau_{\iota,u}=\tilde\tau_{u,\iota}=1_{\tau_u}$ if the initial 
tensor $C^*$--category has conjugates or if every object is a direct sum of irreducibles.\medskip

Let us, informally,  think of $\tau_u\otimes\tau_v$ as a subspace of 
$\tau_{u\otimes v}$.
Equations (A.1) combined with (A.2) require the projection onto
$\tau_{u\otimes v}\otimes\tau_{w}$ to take the subspace 
$\tau_u\otimes\tau_{v\otimes w}$ onto 
$\tau_u\otimes\tau_v\otimes\tau_w$. 
This property should be thought of as a variant of  Popa's {\it commuting square} condition for  
a square of  inclusion of finite
von 
Neumann algebras \cite{Popa2}. In fact in that situation we 
have inclusions $N\subset M$, $Q\subset P$ 
such that $Q\subset N$ and 
$P\subset M$. Recall that this  square is 
called a commuting square 
if $E^M_N(P)\subset Q$ (or, equivalently, if one of the following hold: 
$E^M_P(N)\subset Q$, $E^M_NE^M_P=E^M_PE^M_N=E^M_Q$).\medskip

\noindent
{\bf A.1 Proposition } {\sl Let $(\sigma,\tilde\sigma)$ and $(\tau,\tilde\tau)$ be quasitensor 
functors and suppose $\rho:=\tau\sigma$ is defined. Set $\tilde\rho_{u,v}:=\tau(\tilde\sigma_{u,v})\circ\tilde\tau_{\sigma_u,\sigma_v}$ then $(\rho,\tilde\rho)$ is a quasitensor functor.}\medskip 

\noindent
{\it Proof} The proof just involves routine computations. It is given here for completeness. Obviously $\tilde\rho_{\iota,u}$ and $\tilde\rho_{u,\iota}$ are units for any object $u$.
$$\tilde\rho_{u,v\otimes w}^*\circ\tilde\rho_{u\otimes v,w}=\tilde\tau_{\sigma_u,\sigma_{v\otimes w}}^*\circ\tau(\tilde\sigma_{u,v\otimes w}^*)\circ\tau(\tilde\sigma_{u\otimes v,w})\circ\tilde\tau_{\sigma_{u\otimes v},\sigma_w}=$$
$$\tilde\tau_{\sigma_u,\sigma_{v\otimes w}}^*\circ\tau(1_{\sigma_u}\otimes\tilde\sigma_{v,w})\circ\tau(\tilde\sigma_{u,v}^*\otimes 1_{\sigma_w})\circ\tilde\tau_{\sigma_{u\otimes v},\sigma_w}=$$ 
$$1_{\rho_u}\otimes\tau(\tilde\sigma_{v,w})\circ\tilde\tau_{\sigma_u,\sigma_v\otimes\sigma_w}^*\circ\tilde\tau_{\sigma_u\otimes\sigma_v,\sigma_w}\circ\tau(\tilde\sigma_{u,v}^*)\otimes 1_{\rho_w}=$$
$$1_{\rho_u}\otimes\tau(\tilde\sigma_{v,w})\circ 1_{\rho_u}\otimes\tilde\tau_{\sigma_v,\sigma_w}\circ\tilde\tau_{\sigma_u,\sigma_v}^*\otimes 1_{\rho_w}\circ\tau(\tilde\sigma_{u,v}^*)\otimes 1_{\rho_w}=$$
$$1_{\rho_u}\otimes\tilde\rho_{v,w}\circ\tilde\rho_{u,v}^*\otimes 1_{\rho_w},$$
completing the proof.\medskip 

We comment here on one particularly simple class of quasitensor functors. Let $u$ be an object of a tensor $C^*$-category with irreducible tensor unit, $(\iota,\iota)=\mathbb C$ For an object $u$ pick an orthonormal basis $A_i$ of the Hilbert space $(\iota,u)$ and set $c_u:=\sum_iA_i\circ A_i^*$. $c$ is the support of $\iota$ and is in the centre of $\cA$. Thus, as one sees at once, if $T\in(u,v)$, $T\circ c_u=c_v\circ T$. Note that 
$$c_u\otimes c_v=c_{u\otimes v}\circ 1_u\otimes c_v=c_{u\otimes v}\circ c_u\otimes 1_v.$$ 
In fact 
$$c_{u\otimes v}\circ c_u\otimes 1_v=\sum_ic_{u\otimes v}\circ A_i\otimes 1_v\circ A_i^*\otimes 1_v=\sum_iA_i\otimes c_v\circ A_i^*\otimes 1_v=c_u\otimes c_v.$$

\noindent
{\bf A.2 Proposition } {\sl Let $\tau:\cA\to\cT$ be a $^*$--functor between tensor $C^*$--categories, where $\cT$ has irreducible tensor unit and every object of $\cT$ is a  tensor product of objects in the 
image of $\tau$. Let $u,v$ be objects of $\cA$ and $A_i$ and $B_j$ orthonormal bases of the Hilbert spaces $(\iota,u)$ and $(\iota,v)$, respectively and set 
$$\tilde\tau_{u,v}=\sum_{i,j}\tau(A_i\otimes B_j)\circ\tau(A_i^*)\otimes\tau(B_j^*).$$ 
Then $\tilde\tau$ satisfies all the above conditions except that it may just be a partial isometry. $(\tau,\tilde\tau)$ is a quasitensor functor if and only if $\tau(c_u)=1_{\tau_u}$ for all objects $u$ of $\cA$, $\cT$ is then a full tensor subcategory of a category of Hilbert spaces and $\tilde\tau$ is the unique natural transformation making 
$\tau$ into a quasitensor functor.}\medskip 

\noindent
{\it Proof} It is easily checked that $\tilde\tau$ is a natural transformation and satisfies the associativity condition. Its initial projection is $\tau(c_u)\otimes\tau(c_v)$ and its final projection is $E_{u,v}=\tau(c_u\otimes c_v)$. Hence 
$$E_{u\otimes v,w}\circ E_{u,v\otimes w}=\tau(c_{u\otimes v}\otimes c_w\circ c_u\otimes c_{v\otimes w})=\tau(c_u\otimes c_v\otimes c_w)=E_{u,v,w}.$$ 
Thus $(\tau,\tilde\tau)$ will be a quasitensor functor if and only if $\tau(c_u)=1_{\tau_u}$ for all $u$. In 
particular the support of the tensor unit of $\cT$ is the unit and every object of $\cT$ is a direct sum of 
copies of the unit so $\cT$ is a full tensor subcategory of a category of Hilbert spaces. If $(\tau,\tilde\tau)$ is to be a 
quasitensor functor then $\tilde\tau_{u,v}\circ\tau(A)\otimes\tau(B)=\tau(A\otimes B)$ for all 
$A\in(\iota,u)$ and $B\in(\iota,v)$.\medskip 

The condition $\tau(c_u)=1_u$ for all $u$ is very strong. It also implies that $\tau_u$ is a zero object whenever $c_u=0$, i.e.\ whenever $(\iota,u)=0$. Note that a general quasitensor functor has $E_{u,v}\geq\tau(c_u\otimes c_v)$ with equality characterizing the above special case. 
For this reason, we then say that $\tilde\tau$ is minimal. This case can be alternatively characterized by 
saying that the kernel of $\tau$ is precisely the set of arrows of $\cA$  which are zero when composed with $c$. For if $T\in(u,v)$ and $T\circ c_u=0$ then $\tau(T\circ u)=\tau(T)=0$. Conversely if $\tau(T)=0$ then $\tau(B_j^*\circ T\circ A_i)=0$. Hence $c_v\circ T\circ c_u=T\circ c_u=0$. Thus essentially what the functor $\tau$ does is to map $u$ onto the Hilbert space $(\iota,u)$ and $T$ onto the map $A\mapsto T\circ A$. 

Quasitensor functors $(\tau,\tilde\tau)$ with $\tilde\tau$ minimal are of no direct interest in  this paper as  ${}_\mu\cC_\tau$ reduces to the complex numbers. Indirectly, however, the minimal quasitensor functor $(q,\overline q)$, defined below, plays a role in composition. Let $(\tau,\tilde\tau)$ be a quasitensor functor $\cA\to\cT$ and let $q:\cT\to\cH$ denote the 
$^*$--functor taking an object $x$ of $\cT$ to the Hilbert space $(\iota,x)$ and the arrow $T\in(x,y)$ 
to the map $A\mapsto T\circ A$. There is then a unique quasitensor functor $(q,\tilde q)$ and $\tilde q$ 
is minimal. The composition $(q,\tilde q)\circ(\tau,\tilde\tau)$ is then a quasitensor functor from $\cA$ to $\cH$ without the natural transformation being unitary, in general. When $\tau$ is actually a tensor $^*$--functor, this class of examples was considered  in \cite{PR} and includes as a special case the invariant vectors functor.

Note the following corollaries of the above discussion.\medskip

\noindent
{\bf A.3 Corollary} {\sl Let $\tau:\cA\to\cT$ be a $^*$--functor between tensor $C^*$--categories 
where $\cA$ is a category of Hilbert spaces then $\tau$ may be made into a quasitensor functor in a unique way and is then a relaxed tensor functor 
.}\medskip 

\noindent
{\bf A.4 Corollary} {\sl The tensor product on a $C^*$--category of Hilbert spaces is uniquely defined 
up to a natural unitary transformation. Any two tensor $C^*$--categories of Hilbert spaces are 
equivalent.}\medskip

We next discuss properties of the conjugation on arrows $A\to A^\bullet$ defined in Sect. 4.
This  conjugation does not necessarily commute with the adjoint so that $A^{*\bullet*}$ is, in general, 
an alternative conjugation. If we choose standard solutions of the conjugate equations, however, 
then $A^{*\bullet}=A^{\bullet *}$. In the next lemma, we prove two results that will be used later.\medskip

\noindent
{\bf A.5 Lemma} {\sl Let $A_i$ be an orthonormal basis of the Hilbert space $(\iota,u)$, then $\sum_iA_i\otimes A^\bullet_i=c_u\otimes c_{\overline u}\circ\overline R_u$ and $\sum_iA_i^{*\bullet*}\otimes A_i=
c_{\overline u}\otimes c_u\circ R_u$. Furthermore $c_u^\bullet=c_{\overline u}.$}\medskip 

\noindent
{\it Proof} Let $A,B\in(\iota,u)$ then $A\otimes B^\bullet=A\otimes 1_{\overline u}\circ B^\bullet=(A\circ B^*)\otimes 1_{\overline u}\circ\overline R_u$. Thus $\sum_iA_i\otimes A_i^\bullet=c_u\otimes 1_{\overline u}\circ c_{u\otimes v}\circ\overline R_u=c_u\otimes c_{\overline u}\circ\overline R_u$. The second result can be proved  similarly. Now $c_u^\bullet=\sum_iA_i^\bullet\circ A_i^{*\bullet}.$ Thus  
$c_{\overline u}c_u^\bullet=c_u^\bullet$ and similarly $c_u^\bullet c_{\overline u}=c_{\overline u}$. 
But both $c$ and $c^\bullet$ lie in the centre of $\cA$ and the result follows.\medskip 

A choice of solutions of the conjugate equations determines a scalar product on each $(u,u)$ 
and we write $\phi_u(A^*\circ B):=R^*_u\circ 1_{\overline u}\otimes(A^*\circ B)\circ R_u$. 

If $X\in(\overline u,\tilde u)$ is invertible then $X\otimes 1_u\circ R_u$ and 
$1_u\otimes X^{*-1}\circ\overline R_u$ is another solution of the conjugate equations for $u$. 
Changing $R_u$ using $X$ and $R_v$ using $Y$, the conjugation becomes 
$X^{-1*}\circ A^\bullet\circ Y^*$ whilst the scalar product on $(u,u)$ is given by $\phi_u(A^*\circ B\circ(X^*\circ X)^{\bullet*})$.

Let $R_u,\overline R_u$ and $R_v,\overline R_v$ be solutions of the conjugate equations for $u$ and $v$, then $\overline R_u,R_u$ solves the conjugate equations for $\overline u$ and $1_{\overline v }\otimes R_u\otimes 1_v\circ  R_v,1_u\otimes \overline R_v\otimes 1_{\overline u}\circ\overline R_u$ 
solves the conjugate equations for $u\otimes v$. If these solutions always coincide with $R_{\overline u},\overline R_{\overline u}$ and $R_{u\otimes v},\overline R_{u\otimes v}$ respectively then $u\mapsto R_u$ will be said to be a homomorphic choice of solutions of the conjugate equations.\medskip 

\noindent
{\bf A.6 Proposition} {\sl Let $u\mapsto R_u$ be a homomorphic choice of solutions of the conjugate  
equations then the associated conjugation $^\bullet$ is involutive.}\medskip 

\noindent
{\it Proof} Let $A\in(u,v)$ then $u,v$ are the conjugates of $\overline u,\overline v$ and 
$$A^{\bullet\bullet}=R_{\overline u}^*\otimes 1_v\circ 1_u\otimes(1_{\overline u}\otimes\overline R_v^*\circ 1_{\overline u}\otimes A\otimes 1_{\overline v}\circ R_u\otimes 1_{\overline v})\otimes 1_v\circ 1_u\otimes \overline R_{\overline v}.$$ 
But $R_{\overline u}=\overline R_u$ and $\overline R_{\overline v}=R_v$, so 
$$A^{\bullet\bullet}=\overline R_u^*\otimes 1_v\circ 1_u\otimes(1_{\overline u}\otimes\overline R_v^*\otimes 1_v\circ 1_{\overline u\otimes v}\otimes R_v\circ 1_{\overline u}\otimes A\circ R_u)=$$
$$\overline R_u^*\otimes 1_v\circ 1_{u\otimes u}\otimes A\circ 1_u\otimes R_u=A.$$

The following simple construction shows that, up to an equivalence of tensor $C^*$--categories, we may find such a homomorphic choice. Given a tensor $C^*$--category $\cT$ with a choice $u\mapsto R_u$ of solutions of the conjugate equations, let $\cT^\otimes$ be the tensor $C^*$--category whose objects are words in  the objects $u$ of $\cT$ and their formal adjoints $\overline u$. The tensor product of objects is defined by juxtaposition: $(u_1,u_2,\dots,u_m)\otimes(v_1,v_2,\dots,v_n):=(u_1,u_2,\dots,u_m,v_1,v_2,\dots,v_m)$. The arrows are defined by setting
$$((u_1,u_2,\dots,u_m),(v_1,v_2,\dots,v_n)):=(u_1\otimes u_2\otimes\cdots\otimes u_m,v_1\otimes v_2\otimes\cdots\otimes v_n)$$
and are given the obvious algebraic operations. There is a tensor $^*$--functor $\eta:\cT^\otimes\to\cT$
with $\eta_{(u_1,u_2,\dots u_m)}=u_1\otimes u_2\otimes\cdots\otimes u_m$ an acting as the identity on arrows. Here a formal conjugate $\overline u$ is mapped onto the conjugate of $u$ in $\cT$ determined 
by $R_u$. 
$\eta$ is obviously an equivalence of tensor $C^*$--categories. 
We now choose solutions of the conjugate equations for sequences of 
length one setting $R_{(u)}:=R_u$ and $R_{(\overline u)}:=
\overline R_u$ and extend in the unique way to get a homomorphic choice.\medskip

\noindent
{\bf A.7 Lemma} {\sl Let $u\mapsto R_u$ be a choice of solutions of the conjugate equations for
 $u$ such that
 $$\psi_u(A):={\overline R_u}^*\circ 1_u\otimes A\circ\overline R_u,\quad
  A\in(\overline{u},\overline{u})$$
  is tracial: $\psi_u(A^*\circ B)=\psi_v(B\circ A^*)$, $A,B\in(\overline u,\overline v)$ and 
 write for clarity $J_uS:=S^\bullet$, $S\in(\iota,u)$  then 
$${\rm Tr}(J_u^*\circ J_u)={\rm Tr}(J_u^{-1*}\circ J_u^{-1})={\rm dim}(\iota,u).$$}\medskip 

\noindent
{\it Proof} If $S_i$ is an orthonormal basis of $(\iota,u)$, then by Lemma A.5
$${\rm Tr}(J_u^*J_u)=\sum_{i,j}(S_i\otimes S_i^\bullet)^*\circ(S_j\otimes 
S_j^\bullet)=\overline R_u^*\circ c_u\otimes c_{\overline u}\circ\overline 
R_u=\overline R_u^*\circ 1_u\otimes (c_{\overline u}\circ c_u^\bullet)\circ\overline R_u.$$ 
But by Lemma A.5, $c_u^\bullet=c_{\overline u}$ so we get
 Tr$(J_u^*J_u)=\psi_{ u}(c_{\overline u})$. 
Picking an orthonormal basis $\overline S_i$ for $(\iota,\overline u)$, $c_{\overline u}=\sum_i\overline S_i\circ\overline S_i^*$. But $\psi$ is tracial, so Tr$(J_u^*J_u)=\sum_i\overline S_i^*\circ\overline S_i=$
dim$(\iota,\overline u)=$dim$(\iota,u)$, as required.\medskip

It has been shown in \cite{LR} that if $\psi$ is tracial then the
corresponding conjugation commutes with the adjoint.
As a consequence, if $S,T\in(\iota,u)$ then 
$$(T,S)=(S,T)^*=(S,T)^\bullet=(S^*\circ T)^\bullet=S^{\bullet*}\circ T^\bullet=(S^\bullet,T^\bullet).$$ 
In other words, $^\bullet:(\iota,u)\to(\iota,\overline u)$ is 
antiunitary. \medskip 

Standard solutions $R,\overline R$ of the conjugate equations 
 for $u$ have special properties. 
They are unique up to a unitary and products of standard solutions are again standard. Furthermore, 
$\phi_u$ is independent of the choice of standard solution
so we may replace $\overline R_u$ by $R_{\overline u}$ and get $\psi_u=\phi_{\overline{u}}$ and it is known that $\phi$  is 
tracial, i.e.\ if $A,B\in(u,v)$ then $\phi_u(A^*\circ B)=\phi_v(B\circ 
A^*)$.
Furthermore, $\phi_u(1_u)=R^*\circ R=\overline 
R^*\circ\overline 
R=d(u)$ the intrinsic dimension of $u$
\cite{LR} .\medskip

\noindent
{\bf A.8 Lemma} {\sl Let $u\mapsto R_u$ be a choice of standard solutions of the conjugate equations.  
\begin{description}
\item{a)} Let $W\in(v,u)$ be an isometry then $R_v=W^{\bullet*}\otimes W^*\circ R_u,\overline R_v=W^*\otimes W^{\bullet*}\circ\overline R_u$.
\item{b)} Let $W_i\in(u_i,u)$ be isometries with $\sum_iW_i\circ W_i^*=1_u$, then $R_u=\sum_iW_i^\bullet\otimes W_i\circ R_{u_i}$ and $\overline R_u=\sum_iW_i\otimes W_i^\bullet\circ\overline R_{u_i}$.
\item{c)} There is a unitary $V\in(u,\overline{\overline u})$ such that $R_{\overline u}=V\otimes 1_{\overline u}\circ \overline R_u$ and $\overline R_{\overline u}=1_{\overline u}\otimes V\circ R_u$.
\end{description}}\medskip

\noindent 
{\it Proof}\begin{description} 
\item a) is proved as follows
$$\overline R_v=W^*\otimes 1_{\overline v}\circ 1_u\otimes W^{\bullet*}\circ\overline R_u=
W^*\circ W\otimes 1_{\overline v}\circ\overline R_v=R_v.$$ 
The second equation can be proved similarly.
\item To prove b), we compute as follows. 
$$\sum_iW_i^\bullet\otimes W_i\circ R_{u_i}=\sum_i1_{\overline u}\otimes (W_i\circ W_i^*)\circ R_u=R_u$$ 
and the second equation can be proved similarly.
\item To prove c) note that $\overline R_u,R_u$ is a standard solution of the conjugate equations for $\overline u$ and therefore differs from $R_{\overline u},\overline R_{\overline u}$ by a unitary $V$ as 
claimed.
\end{description}\medskip

 It must be remembered though that even if $u\to R_u$  is standard, $\hat R_u$ and $\tilde R_u$  defined in Sect. 4 will not, in general, be standard, nor will the conjugation commute with the adjoint. Nevertheless the following result holds.\medskip 

\noindent
{\bf A.9  Lemma} {\sl Let $u\mapsto R_u$ be a standard choice of solutions of the conjugate equations 
then the set of objects $u$ such that $\hat R_u,\hat{\overline R}_u$ is a standard solution of the 
conjugate equations for $\mu_u$ is closed under tensor products subobjects, direct sums and conjugates.}\medskip

\noindent
{\it Proof} Suppose $\hat R_u,\hat{\overline R}_u$ is standard and $W\in(v,u)$ is an isometry then 
by Lemma A.8, $\hat R_v=\mu(W^{\bullet*})\otimes\mu(W^*)\circ\hat R_u$.
Now 
$$\hat R_v^*\circ \hat R_v=\hat R_u^*\circ\mu(E^\bullet)\otimes\mu(E)\circ\hat R_u=\hat R_u^*\circ 1_{\mu_{\overline u}}\otimes\mu(E)\circ\hat{\overline R}_u=\phi_{\mu_u}(E),$$ 
where $\phi_u$ is the standard left inverse of $\mu_u$. By the tracial property of the standard left inverse, $\phi_{\mu_u}(E)=\phi_{\mu_u}(W^*\circ W)=\phi_{\mu_v}(1_{\mu_v})=d(\mu_v)$. Similarly 
$\hat{\overline R}_v^*\circ \hat{\overline  R}_v=d(\mu_v)$  and $\hat R_v,\hat{\overline R}_v$ are standard. Now suppose $\hat R_{u_i},\hat{\overline R}_{u_i}$ are standard and $W_i\in(u_i,u)$ are 
isometries with $\sum_iW_i\circ W_i^*=1_u$ then by Lemma A.7, $\hat R_u=\sum_iW_i^\bullet\otimes W_i\circ\hat R_{u_i}$ hence $\hat R_u^*\circ\hat R_u=\sum_i\hat R_{u_i}^*\circ\hat R_{u_i}=\sum_id(\mu_{u_i})=d(\mu_u)$. Similarly, $\hat{\overline R}_u^*\circ\hat{\overline R}_u=d(\mu_u)$ so that $\hat R_u,\hat{\overline R}_u$ is standard. Again, if $\hat R_u$ is standard, by Lemma A.7, there is a unitary $V\in(u,\overline{\overline u})$ such that $\hat R_{\overline u}=\mu(V)\otimes 1_{\mu_{\overline u}}\circ\hat{\overline R}_u$ thus $\hat R_{\overline u}^*\circ 
\hat R_{\overline u}=\hat{\overline R}_u^*\circ\hat{\overline R}_u=d(\mu_u)=d(\mu_{\overline u})$ and similarly $\hat{\overline R}_{\overline u}^*\circ\hat{\overline R}_{\overline u}=d(\mu_{\overline u})$. Thus $\hat R_{\overline u}$ is standard. 
The question of whether $\hat R_u$ is standard is obviously independent of the choice of standard solution $R_u$. If $R_{u\otimes v}$ is chosen to be of product form then 
the same is true of $\hat R_{u\otimes v}$. Thus $\hat R_u$ and $\hat R_v$ standard 
imply that $\hat R_{u\otimes v}$ is standard.\medskip 

\noindent
{\bf A.10 Corollary} {\sl If $u\mapsto R_u$ is standard and $\mu_v$ is an irreducible generator of ${\cal M}$ 
then $u\mapsto\hat R_u$ is standard. If $d(u)=d(\mu_u)$ then $\hat R_u,\hat{\overline R}_u$ is 
standard.} 

\noindent
{\it Proof} The first statement follows since, $\mu_v$ being irreducible, $\hat R_v,\hat{\overline R}_v$ is 
automatically standard. Now $d(\mu_u)\leq\hat R_u^*\circ \hat R_u=\mu(R_u^*)\circ E_{\overline u,u}\circ\mu(R_u)\leq R_u^*\circ R_u=d(u)=d(\mu_u)$ and the result follows.\medskip

\end{section}

\noindent
{\it Acknowledgement} We gratefully acknowledge the numerous 
discussions with S. Doplicher on the contents of this paper.


\begin{thebibliography}{VD}

\bibitem{Banica0} T. Banica: Th\'eorie des repr\'esentations du groupe 
quantique compact libre $O(n)$, {\it C.R. Acad. Sci. Paris}, {\bf 322}
(1996), 241--244.

\bibitem{Banica} T. Banica: Le groupe quantique compact libre $U(n)$,
{\it Comm. Math. Phys.}, {\bf 190} (1997), 143--172.

\bibitem{BRV} J. Bichon, A. De Rijdt, S. Vaes: Ergodic coactions with
large multiplicity and monoidal equivalence of quantum groups, {\it Comm.
Math. Phys.}, {\bf 262} (2006), 703--728.

\bibitem{Boca} F. Boca: Ergodic actions of compact matrix pseudgroups  on
$C^*$--algebras, {\it Ast\' erisque}, {\bf 232} (1995), 93--109.



\bibitem{DRInventiones} S. Doplicher, J.E. Roberts: A new duality theory
for compact groups, {\it Invent. Math.}, {\bf 98} (1989),
157--218.

 
\bibitem{GW} F.M. Goodman, H. Wenzl: Ideals in the Temperley-Lieb category. 
Appendix to M. Freedman: a magnetic model with a possible Chern-Simons phase. 
{\it Comm.\ Math.\ Phys.\ }, {\bf 234}, (2003), 129-183.


\bibitem{Haag} R. Haag: Local quantum physics. Fields, particles, 
algebras. Second edition. Texts and Monographs in Physics. 
Springer-Verlag, Berlin, 1996. 

\bibitem{HLS} R. H\o egh--Krohn, M. Landstad, E.
St\o rmer: Compact ergodic groups of automorphisms, {\it Ann. of Math.},
{\bf 114} (1981) 137--149.

 \bibitem{LR} R. Longo, J.E. Roberts: A theory of dimension, {\it
$K$--Theory}, {\bf 11} (1997), 103--159.

 
\bibitem{PR} C. Pinzari, J.E. Roberts: A duality theorem for ergodic 
actions of compact quantum groups on $C^*$--algebras, {\it Comm. Math. Phys.}, 
{\bf 277}, (2008). 385-421.

\bibitem{PR1} C. Pinzari, J.E. Roberts: Ergodic actions of $S_\mu U(2)$ on $C^*$--algebras from $II_1$ subfactors, arXiv:math0806.4519v1

\bibitem{PRinduction} C. Pinzari, J.E. Roberts: A theory of induction and classification of tensor $C^*$--categories,
{\it J. Noncomm. Geom.} {\bf  6}, (2012),  665--719.

\bibitem{PRlinking} C. Pinzari, J.E. Roberts: Linking $C^*$--bimodules arising from pairs of functors, in progress.


\bibitem{Podles} P. Podles: Symmetries of quantum spaces. Subgroups and
quotient spaces of quantum $SU(2)$ and $SO(3)$ groups, {\it Comm. Math.
Phys.}, {\bf 170} (1995), 1--20.

 
\bibitem{Popa2} S. Popa: Maximal injective subalgebras in factors 
associated with free groups, {\it Adv. Math.}, {\bf 50} (1983), 27--48.

\bibitem{Tomatsu} R. Tomatsu: A characterization of right coideals of
quotient type and its application to classification of Poisson boundaries,
{\it Comm.\ Math.\ Phys.\ } {\bf 275}, (2007), 271-296.


\bibitem{WVD} A. Van Daele, S. Wang: Universal quantum groups, {\it 
Internat. J. Math.}, {\bf 7} (1996), 255--263.

\bibitem{Wang}  S. Wang: Ergodic actions of universal quantum groups on
operator algebras, {\it Comm. Math. Phys.}, {\bf 203} (1999), 481--498.

\bibitem{Wang2} S. Wang: Structure and isomorphism classification of 
compact quantum groups $A_u(Q)$ and $B_u(Q)$, {\it J. Operator Theory},
{\bf 48} (2002), 573--583.

\bibitem{Wassermann1} A. Wassermann: Ergodic actions of compact groups
on operator algebras. I. General theory. {\it Ann. of Math.}, {\bf 130} 
(1989),
273--319.


\bibitem{Wassermann2} A. Wassermann: Ergodic actions of compact groups
on operator algebras. II. Classification of full multiplicity ergodic
actions. {\it Canad. J. Math.}, {\bf 40} (1988),
1482--1527.


\bibitem{Wassermann3} A. Wassermann: Ergodic actions of compact groups
on operator algebras. III. Classification for $SU(2)$. {\it Invent.  
Math.}, {\bf 93} (1988),
309--354.


 \bibitem{Wenzl3} H. Wenzl: $C^*$--tensor categories from quantum groups,
{\it J. Amer. Math. Soc.}, {\bf 11} (1998), 261--282.

\bibitem{Wcmp} S. L. Woronowicz: Compact matrix pseudogroups, {\it Comm.
Math. Phys.}, {\bf 111} (1987), 613--665.

\bibitem{WoronowiczTK} S. L. Woronowicz: Tannaka--Krein duality for 
compact matrix
pseudogroups. Twisted $SU(N)$ groups, {\it Invent. Math.}, {\bf 93}
(1988), 35--76.

\bibitem{WLesHouches} S. L. Woronowicz: Compact quantum groups, Les 
Houches, 1995,
845--884, North Holland, Amsterd
1998.

\bibitem{Y} S. Yamagami: A categorical and diagrammatical approach to Temperley-Lieb algebras, 
arXiv:math/0405267v2.

\bibitem{Y1} S. Yamagami: Oriented Kauffman diagrams and universal quantum groups, 
arXiv:math/0602117v2.

\bibitem{Y2} S. Yamagami: Fiber functors on Temperley-Lieb categories, arXiv:math/0405517.




\end{thebibliography}
\end{document}